\documentclass[11pt,psamsfonts,twoside]{article}
\usepackage{mathrsfs}
\usepackage{bbm}
\usepackage{amsmath}
\usepackage{amsmath,amsthm}
\usepackage{amssymb,amscd}
\usepackage{amsfonts,amsbsy}
\usepackage{fancyhdr,graphicx}
\usepackage[dvips]{psfrag}
\usepackage{indentfirst}
\usepackage{xcolor}
\usepackage{setspace}

\numberwithin{equation}{section}
\newtheorem {proposition} {Proposition}[section]
\newtheorem {theorem}     [proposition]{Theorem}
\newtheorem {corollary}   [proposition]{Corollary}
\newtheorem {lemma}       [proposition]{Lemma}

\newtheorem {remark}      [proposition]{Remark}

\newtheorem {definition}  [proposition]{Definition}
\allowdisplaybreaks

\newcommand{\supp}{{\rm supp}}
\newcommand{\diva}{{\rm div}}

\begin{document}
\setlength{\parindent}{4ex} \setlength{\parskip}{1ex}
\setlength{\oddsidemargin}{12mm} \setlength{\evensidemargin}{9mm}
\title{Well-posedness and ill-posedness of the 3D generalized Navier-Stokes equations in Triebel-Lizorkin spaces}
\author{Chao Deng,\quad \quad Xiaohua Yao}

\date{}
\maketitle
\begin{abstract}
In this paper, we study the Cauchy problem of the 3-dimensional (3D)
generalized incompressible Navier-Stokes
 equations (gNS) in Triebel-Lizorkin space $\dot{F}^{-\alpha,r}_{q_\alpha}(\mathbb{R}^3)$ with
 $(\alpha,r)\in(1,\frac{5}{4})\times[2,\infty]$ and $q_\alpha=\frac{3}{\alpha-1}$. Our work establishes a {\it dichotomy} of well-posedness and ill-posedness
 depending on $r=2$ or $r>2$. Specifically, by combining the new endpoint bilinear estimates in $L^{\!q_\alpha}_x\!L^2_T$ with the characterization of
 Triebel-Lizorkin space via fractional semigroup, we prove the  well-posedness of the gNS in $\dot{F}^{-\alpha,r}_{q_\alpha}(\mathbb{R}^3)$ for $r=2$.
 On the other hand, for any $r>2$, we show that the solution to the gNS can develop {\it norm inflation} in the sense that arbitrarily small initial data
 in the spaces $\dot{F}^{-\alpha,r}_{q_\alpha}(\mathbb{R}^3)$ can lead the corresponding solution to become arbitrarily large after an arbitrarily short time.
 In particular, such dichotomy of Triebel-Lizorkin spaces is also true for the classical N-S equations, i.e.\,\,$\alpha\!=1$.
Thus the Triebel-Lizorkin space framework  naturally provides better
 connection between the well-known Koch-Tataru's $BMO^{-1}$
 well-posed work and Bourgain-Pavlovi\'c's $\dot{B}_\infty^{-1,\infty}$ ill-posed work.
 \end{abstract}
 {\small \noindent{\bf Keywords:}\, Generalized Navier-Stokes equations; Triebel-Lizorkin space; well-posedness; Ill-posedness.}
\maketitle
\section{Introduction}
 In this article, we study the initial value problem of the following 3D generalized incompressible Navier-Stokes equations (gNS):
   \begin{align}\label{DLY-1}
  \left\{\begin{aligned}
   &\partial_t u + (-\Delta)^{\alpha} u + u\cdot \nabla u +\nabla
   p=0,\\
   &\nabla\cdot u=0,\\
   &u(x,0)=u_0(x),
     \end{aligned}\right.
 \end{align}
 where $\alpha >{1\over 2} $, $(x,t)\in\mathbb{R}^3\times(0,\infty)$, ${u}(x,t)=\!(u^1(x,t), u^2(x,t),\!u^3(x,t))$
 are unknown vector functions, $p(x,t)$ is unknown scaler function,  and ${u}_0(x)$ is a given vector function satisfying
 divergence free condition $\nabla\cdot{u}_0=0$.

 Mathematical analysis of the classical incompressible Navier-Stokes (N-S) equations (\,$\alpha=1$) has a long history. It goes back to
 Leray's famous work, i.e. \cite{LER34}, in which Leray first introduced the concept of weak solutions and proved existence
 of global weak solutions associated with $L^2(\mathbb{R}^n)$ initial data by using an approximation approach
 and some weak compactness arguments. In 1964, Fujita-Kato \cite{FUJK64} initiated a different approach and proved well-posedness of the
 initial value problem of the N-S in ${H}^s(\mathbb{R}^n)$ for $s\geq\frac{n}{2}-1$. This approach was later extended to
 various other function spaces, see \cite{CAN95,CAN04,LEM02,PLA96} for expositions and references therein. An
interesting result that must be mentioned is due to Koch-Tataru
\cite{KOCT01}. They proved that the solutions of N-S are well-posed
in $BMO^{-1}$ which is the largest function space for
well-posedness.
 Besides these well-posedness results in critical spaces mentioned above, there also exist several works for supper-critical initial
 value, for instance, existence of
 solutions for the initial value problem of the N-S for initial
 data in supercritical spaces $L^2(\mathbb{R}^n)$ space (which is supercritical for $n\geq 3$ but
 critical for $n=2$) and sums of $L^2(\mathbb{R}^n)$ with some well-posedness spaces (cf. \cite{CAL90,LEM02}).

 As we know, one crucial reason of working with the gNS equations on $\mathbb{R}^3$ for $\alpha>\frac{1}{2}$ is that they
 provides us the deeper understandings of the different actions of
 fractional Laplacian.  Similar to the classical Navier-Stokes equations, one of the most
  primary problems is to establish local or global-in-time well-posedness of the gNS equation. Do the solutions exist in some spaces?  If so, are they unique and
  is the system stable for certain initial data?  By stable we mean
  that small perturbation of initial data guarantees small perturbation of solution. It is worth mentioning that either stability or instability
  of the nonlinear PDEs
  has a lot of applications in numerical analysis field.

 Up to now,  there exist many interesting works about the
 well-posedness and ill-posedness of gNS equations.
 In Table 1.1, we list several important progresses in the Besov
space framework:
 \vspace*{-1ex}
 \begin{spacing}{1.4}\begin{center}
 \begin{tabular}{|r|l|}
 \hline
\multicolumn{2}{|c|}
{Well-posedness/\,ill-posedness for the 3D generalized Navier-Stokes equations}  \\
   \hline
$\frac{1}{2}<\alpha<1$  &   Well-posed in $\dot{B}^{1-2\alpha,\infty}_\infty$\  \cite{LZ10, YuZhai:2012}.   \\
\hline
           $\alpha=1$           &   Well-posed in {\small $BMO^{-1}$}\ \cite{KOCT01},  Ill-posed in $\dot{B}^{-1,\infty}_\infty$\  \cite{BP08} and                             \\
                                &   ill-posed in the logarithm-type Besov space \  \cite{Yoneda:2010}.
                                \\ \hline
 $1<\alpha<\frac{5}{4}$ &   Well-posed in { $\dot{B}^{s,\,r}_{2}$}\! with $s=1\!-\!2\alpha+\frac{3}{2}$ and $1\le r\!\le\infty$\ \cite{Wu:2006803},\\
                                &  Ill-posed in $\dot{B}^{1-2\alpha,\infty}_\infty$\  \cite{CheskidovaShvydkoyb:2012}. \\
 \hline $\frac{5}{4}\le \alpha$     &   Global existence and uniqueness of classical solutions\  \cite{Lions:1969}.\\
 \hline
 \end{tabular}\\
\vskip0.3cm
 { Table 1.1}
\end{center}
\vspace*{-3ex}
\end{spacing}
 Another result must be mentioned is that Cheskidov and Dai \cite{cheskidovdai:2012} proved norm-inflation of the gNS in subcritical
 Besov spaces $\dot{B}^{-\alpha,\infty}_\infty$ for $\alpha>1$.

All the function spaces in Tabel 1.1 are {\it invariant under
 scaling}
 $
 {u_{0\lambda}}\!=\!\lambda^{2\alpha\!-1}u_0(\lambda x)
 $,
which corresponds to the solution $u$ of the gNS scale invariant
under transformation
\begin{align}\label{scale}
{u}_\lambda(x,t)\!=\!\lambda^{2\alpha-1}{u}(\lambda
x,\lambda^{2\alpha}t) \;\;\text{ with }
 \lambda>0 \text{ and } \alpha>{1}/{2}.\end{align}
 By applying \cite[Proposition 4.2]{LEM02} to \eqref{scale}, we can check that $\dot{B}^{1-2\alpha,\infty}_\infty$  is the largest scale invariant
 function space for the gNS equations.  From Table 1.1, we observe that,
 when ${\frac{1}{2}}<\!\alpha<\!1$, the gNS equations are  well-posed in the largest scale invariant
 function space $\dot{B}^{1-2\alpha,\infty}_\infty$ (\,\cite{LZ10,YuZhai:2012}\,), but when $1\le \!\alpha<\!\frac{5}{4}$, the gNS equations actually
 are  ill-posed in $\dot{B}^{1-2\alpha,\infty}_\infty$ (\,\cite{BP08, CheskidovaShvydkoyb:2012}). Hence for the cases  $1\le \!\alpha<\!\frac{5}{4}$,
 there exists certain difference between well-posed space with the largest scale invariant space $\dot{B}^{1-2\alpha,\infty}_\infty$.

 In fact, when $\alpha=1$ ( i.e. the  classical N-S equation ), the best well-posed space is
  $\dot{F}^{-1,2}_{\infty}=BMO^{-1}$ given by Koch-Tataru \cite{KOCT01}. Furthermore,  $\dot{F}^{-1,2}_{\infty}\varsubsetneq
 \dot{B}^{-1,\infty}_\infty$ and there does exist a minor difference between $\dot{F}^{-1,2}_{\infty}$ and
 $\dot{B}^{-1,\infty}_{\infty}$.
 When $1<\!\alpha< \frac{5}{4}$, the known results on well-posedness, for instance, $B^{s,r}_2$ with $1\le r\le\infty$ and $s={1-2\alpha+}\frac{3}{2}$ (\,see\,\cite{Wu:2006803}), are surely
not the largest well-posed space of the gNS  from \eqref{B-3} (\,see
Subsection 2.3 below\,). Hence it would be very interesting to
figure out  what is the largest well-posed space and examine how
large the difference is with  $\dot{B}^{1-2\alpha,\infty}_\infty$.
Roughly speaking, we guess the difference between the largest
well-posed space and the
     largest scale invariant space $\dot{B}^{1-2\alpha,\infty}_\infty$ should enlarge as $\alpha$ increases.

 In this paper, we investigate the interesting problem in  critical Triebel-Lizorkin space framework  and  obtain
 the following conclusions about the 3D gNS with
 $1\le\alpha<\frac{5}{4}$ in the whole space $\mathbb{R}^3$:
 \begin{spacing}{1.5}
 \vspace*{-1ex}
\begin{center}
 \begin{tabular}{|r|l|}
 \hline
\multicolumn{2}{|c|}
{Well-posedness and ill-posedness of 3D gNS in critical Triebel-Lizorkin spaces}  \\
   \hline
  $\alpha=1$   & well-posed in $\dot{F}^{-1,2}_{\infty}$ (\,\cite{KOCT01}\,) and ill-posed in $\dot{F}^{-1,r>2}_{\infty}$ (Remark 1.4) \\
 \hline
 $1\!<\!\alpha\!<\!\frac{5}{4}$&  well-posed in $\dot{F}^{-\alpha,2}_{q_\alpha}$ and ill-posed in $\dot{F}^{-\alpha,r>2}_{q_\alpha}$
 (Theorems 1.2--1.3\,) \\
   \hline
 \end{tabular}\\
\vskip0.3cm
  { Table 1.2}
\end{center}
\vspace*{-3ex}
\end{spacing}

  Table 1.2 shows that  the largest Triebel-Lizorkin-type well-posed space for gNS is
  $\dot{F}^{-\alpha,2}_{q_\alpha}$ for any $1\le\alpha<\frac{5}{4}$.  Indeed,  for any $r<2$,
  we have $ \dot{F}^{-\alpha,r}_{q_\alpha}\!\subsetneq\dot{F}^{-\alpha,2}_{q_\alpha}$ and  for any $r>2$,
  we prove that the gNS is ill-posed in $\dot{F}_{q_\alpha}^{-\alpha,r}(\mathbb{R}^3)$
  in the sense that arbitrarily small initial data can lead the corresponding solution to become arbitrarily large after an arbitrarily short time.
 Hence  the  results in Table 1.2  sharpen the well-posed
 analysis of gNS in  the critical Triebel-Lizorkin space
 $\dot{F}^{-\alpha, 2}_{q_\alpha}$ for all $\alpha\in[ 1, {5\over
 4})$.

 The results obtained here and \cite{cheskidovdai:2012} indicate that ill-posedness is more closely related to the smoothing indices $s$ than
 the integrability indices $q$. More precisely, the gNS is ill-posed in $\dot{B}^{-\alpha,\infty}_{q}$ for $q_\alpha\le q\le\infty$.

 Next  we would like to introduce the {\it ideas} of the paper.
 In most applications, people use space-time type norm where one
 takes the space norm first. Here, we use some time-space type norm.
 This seems to be critical since it seems impossible to get the best
 results without these new norms and the related estimates in $L^{q_\alpha}_xL^2_T$ in which we have an equivalent characterization of
 Triebel-Lizorkin space by fractional semigroup (see Appendix A
 below). Then combining the a-priori bilinear estimates and the
 contraction arguments we prove the well-posedness.  To show the ill-posedness, we shall adopt the novel
 framework of norm inflation first introduced by
 Bourgain-Pavlovi{\'c} \cite{BP08} in their study of the ill-posedness of the Navier--Stokes equation in
 $\dot{B}_{\infty}^{-1,\infty}(\mathbb{T}^3)$; but in doing so, we introduce some new inputs to the gNS. In particular,
  we make use of Hardy-Littlewood maximal function to estimate the norm of the solution in Triebel-Lizorkin
 space.
 Finally, we  {\it conjecture} that for the specially constructed initial data in Subsection 3.2 below, there should exist
 a unique global classical solution since the date is not only energy finite but also essentially 2-dimensional/ summation of plane waves.

 In order to prove the main results, we first recall the definition of homogeneous Besov
 spaces and Triebel-Lizorkin spaces (\cite{Triebel:1978}).
 Let $\varphi(\xi)=\varphi(|\xi|)$ be a real-valued  smooth function such that
 $0\le \varphi(\xi)\le1$ and
  \begin{align}\label{DLY-1.4}
  \supp\varphi\subset\{\xi\in\mathbb{R}^3;\; {3}/{4}\le|\xi|\le  {8}/{3}\},\;\;\;\text{
                  } \sum_{j\in\mathbb{Z}}\varphi(2^{-j}\xi)=1\
               \text{ for any } \xi\ne0.
  \end{align}
 For any tempered distribution $f$ and $i,j\in\mathbb{Z}$, we define the dyadic block as follows:
 \begin{align}\label{DLY-1.5}
 \Delta_{j}f(x)=\varphi(2^{-j}\nabla)f(x)\quad \text{ and\;\; }
  \Delta_{i}\Delta_jf\equiv0 \text{ if } |i-j|\ge2.
\end{align}

 To exclude nonzero polynomials in homogeneous Besov spaces and Triebel-Lizorkin spaces,
 it is natural to use $Z'(\mathbb{R}^3)$ to denote the subspace of tempered
 distribution $f\in \mathcal{S}'(\mathbb{R}^3)$ modulo all polynomials set $P(\mathbb{R}^3)$, i.e.
 $$Z'(\mathbb{R}^3)=S'(\mathbb{R}^3)/P(\mathbb{R}^3).$$

 Now we give the definition of Besov and Triebel-Lizorkin spaces, see
\cite{Triebel:1978}.
\begin{definition}\label{Def:1.1}
 For any $(s,q,r)\in(-\infty,\infty)\times[1,\infty]\times[1,\infty]$, we denote by $\dot{B}^{s,r}_{q}(\mathbb{R}^3)$ the set of distribution
  $f\in Z'(\mathbb{R}^3)$ satisfying
 \begin{align*}
   \|f\|_{\dot{B}^{s,r}_{q}(\mathbb{R}^3)}=\|\{2^{js}\|\Delta_jf\|_{L^q_x}\}\|_{l^r}<\infty
 \end{align*}
 and we let $\dot{F}^{s,r}_{q}(\mathbb{R}^3)$ with $1\le q<\infty$ be the set of distribution $f\!\in\!Z'(\mathbb{R}^3)$ satisfying
 \begin{align*}
   \|f\|_{\dot{F}^{s,r}_{q}(\mathbb{R}^3)}=\|\|\{2^{js}\Delta_jf\}\|_{l^r}\|_{{L^q_x}}<\infty.
 \end{align*}
 In particular, $\dot{F}^{s,r}_\infty$
  is defined by the usual modification as in {\rm
\cite[Definition, p.30]{Triebel:1978}}.
\end{definition}

We are ready to state our main results on well-posedness and ill-posedness. 

 \begin{theorem}\label{Thm:1.2}\!$(\mathrm{Well\text{-}posedness})$
   Let $1<\alpha<\frac{5}{4}$ and
   $q_\alpha=\frac{3}{\alpha-1}$. Then for any
   $u_0\in\!\dot{F}^{-\alpha,2}_{q_\alpha}(\mathbb{R}^3)$ satisfying
   $\nabla\cdot u_0=0$,
   there exists $T\!=\!T(u_0)\!>\!0$ such that system \eqref{DLY-1} has a unique local solution $u$ satisfying
  \begin{align*}
     u \in C([0,T];\dot{F}^{-\alpha,2}_{q_\alpha}(\mathbb{R}^3))\cap
     L^{q_\alpha}_xL^2_T.
  \end{align*}
  Furthermore, if $\|u_0\|_{\dot{F}^{-\alpha,2}_{q_\alpha}(\mathbb{R}^3)}$ is small enough,
  then system \eqref{DLY-1} has a unique global solution satisfying
  \begin{align*}
    u \in C([0,\infty);\dot{F}^{-\alpha,2}_{q_\alpha}(\mathbb{R}^3))\cap
    L^{q_\alpha}_xL^2_t.
  \end{align*}
 \end{theorem}

 \vskip0.4cm

\begin{theorem}\label{Thm:1.3}\!$(\mathrm{Ill\text{-}posedness})$
 For any $1<\alpha<\frac{5}{4}$, $r>2$, $\delta>0$ and $q_\alpha=\frac{3}{\alpha-1}$,  there exists a solution $u$ to  system \eqref{DLY-1} with initial data $u_0\in
\dot{F}^{-\alpha,r}_{q_\alpha}(\mathbb{R}^3)$ satisfying
\begin{align*}
 \|u_0\|_{\dot{F}^{-\alpha,r}_{q_\alpha}(\mathbb{R}^3)}\lesssim\delta
  \end{align*}
 and $ \nabla\cdot u_0=0$ such that for some $0<T<\delta$,
\begin{align*}
  \|u(T)\|_{\dot{F}^{-\alpha,r}_{q_\alpha}(\mathbb{R}^3)}\gtrsim\frac{1}{\delta}.
\end{align*}
\end{theorem}

\smallskip
\begin{remark}
\label{Rem:1.4}
 {\rm \hskip12cm \;\;

  (i)\, In the proof of Theorem \ref{Thm:1.3}, the constructed
 initial datum satisfy several {\it good} properties: {\it real-valued,
 smooth, energy finite in the whole space, and essentially plane
 waves (almost 2D)}.\;

 (ii) When $\alpha=1$ and $q_\alpha=\infty$,
$\dot{F}^{-1,2}_\infty=BMO^{-1}$,  by  the { Koch-Tataru's
well-posed work \cite{KOCT01} in
 $BMO^{-1}$},  we also  proved the {\it ill-posedness } of the Navier-Stokes equations in the
   Triebel-Lizorkin spaces $\dot{F}^{-1,r}_\infty(\mathbb{R}^3)$ (\,cf.\,\cite{DengYao:2012}\,)  for  $2<r<\infty$ which are strictly smaller than $\dot{B}^{-1,\infty}_\infty(\mathbb{R}^3)$  in which
    Bourgain and Pavlovi\'c proved ill-posedness of the Navier-Stokes equations, see
    \cite{BP08}.

 }
 \end{remark}

\smallskip

 { This paper is {\it organized} as follows:}
 In Section 2, we mainly establish the well-posedness by proving
 a key bilinear estimate on $L^{q_\alpha}_xL^2_T$.  Moreover, we also consider many other bilinear estimates in the end of this section. In particular, by interpolation several applications of our bilinear estimate are given; In Section 3, we
 first construct a very special initial data and list some
 necessary remarks, and then we establish all the desired estimates
 about the first and second approximation terms which will be used in controlling the remainder term.
 Finally, combining all the a-priori estimates we prove ill-posedness of the gNS.

  \smallskip
 \noindent\textit{Notations}: Throughout this paper, we shall use $C$ and $c$ to denote universal constants
  and may change from line to line. Both $\mathcal{F}f$ and $\widehat{f}$  stand for Fourier
 transform of $f$ with respect to space variable, while
 $\mathcal{F}^{-1}$ stands for the inverse Fourier transform.  We
 denote $A\le{CB}$ by $A\lesssim B$ and $A\lesssim{B}\lesssim{A}$
  by $A\sim{B}$. For any $1\le p\le\infty$, we denote
   $L^p(0,T)$, $L^p(T_1,T_2)$, $L^p(0,\infty)$ and $L^q(\mathbb{R}^3)$ by
   $L^p_T$, $L^p_{[T_1,T_2]}$, $L^p_t$ and $L^q_x$, respectively.
 Later on, we use $\dot{F}^{s,r}_q$ to denote $\dot{F}^{s,r}_q(\mathbb{R}^3)$
 if there is no confusion about the domain, and similar conventions
 are applied.

\section{Analysis of well-posedness}
 In this section, we will prove well-posedness of the 3D gNS in $\dot{F}^{-\alpha,2}_{q_\alpha}$
 for
 $q_\alpha=\frac{3}{\alpha-1}$ and $1<\alpha<\frac{5}{4}$.
 Notice that for $\alpha=1$, $\dot{F}^{-\alpha,2}_{q_\alpha}=BMO^{-1}$
 in which the well-posedness is proved by Koch and Tataru, see \cite{KOCT01}.

 As usual, we first write \eqref{DLY-1} into
 the following equivalent mild integral equations:
\begin{align}
  &u=e^{-t(-\Delta)^\alpha}u_0-\int_{0}^te^{-(t-\tau)(-\Delta)^\alpha}\mathbb{P}(u\cdot\nabla v) d\tau,\label{DLY-2}
  \end{align}
  where $\mathbb{P}$ is the Leray projection operator and
  $\mathbb{P}=Id-\nabla\frac{1}{\Delta}\diva$.

 For simplicity, we denote the bilinear term by
\begin{align}
  \label{DLY-3}
 B(u,v):&=\int_{0}^te^{-(t-\tau)(-\Delta)^\alpha}\mathbb{P}(u\cdot\nabla v) d\tau.
 \end{align}

 To prove well-posedness, we first prove several preliminary lemmas including the endpoint bilinear estimate.
 Based on these estimates,  the well-posedness immediately follows from  the
 standard Picard iteration principle.
 Finally,  we also give some other bilinear estimates and  their applications.

 \subsection{Preliminaries}
  In this subsection, we first give several preliminary lemmas.
  The first lemma is about the point-wise estimates for the kernel of fractional semigroup $e^{-t(-\Delta)^\alpha}$ with a regularized operator
  $(-\Delta)^{s/2}$. For convenience of other applications, we consider dimension $n\ge 2$ and also allow $\alpha>0$.  For any $s>-n$, $j\in \mathbb{Z}$, $t>0$
  and $x\in\mathbb{R}^n$, let
$$K_s^\alpha(t,x)=\mathcal{F}^{-1}\Big(|\xi|^s\ e^{-t|\xi|^{2\alpha}}\Big)(x), $$
and
$$\hskip1.4cm K_{s,j}^\alpha(t,x)=\mathcal{F}^{-1}\Big(|\xi|^s\ e^{-t|\xi|^{2\alpha}}\varphi(2^{-j}\xi)\Big)(x).  $$
where $\varphi(\xi)$ be a smoothing truncation function defined in
\eqref{DLY-1.4}.  Clearly, they are the kernel of the operator
families $(-\Delta)^{s/2} e^{-t(-\Delta)^\alpha}$ and
$(-\Delta)^{s/2} e^{-t(-\Delta)^\alpha}\Delta_j$, respectively.  As
$s=0$, we denote the kernels by $K^\alpha(t,x)$ and
$K^\alpha_j(t,x)$, respectively. \vskip0.3cm
\begin{lemma}\label{Lem:2.2} Let $\alpha>0$, $n\ge2$, $s>-n$ and $\beta\in\mathbb{N}^n$. Then we have the following estimates:
\begin{align}\label{pointwiseA}
 |\partial^\beta_x K^\alpha_s(t,x)|\le
 \left\{
 \begin{aligned}
 & C\ (t^{\frac{1}{2\alpha}}+|x|)^{-(n+s+|\beta|)},  &s\neq0;
 \\
 &C\ t (t^{\frac{1}{2\alpha}}+|x|)^{-(n+2|\alpha|+|\beta|)},  \ &s=0.
 \end{aligned}\right.
 \end{align}
and for any $N\ge1$,
\begin{align}\label{pointwiseB}
   |\partial^\beta_x K^\alpha_{s,j}(t,x)|\le C_N\  e^{-c t2^{2j\alpha}}2^{j(n+s+|\beta|)}\ (1+|2^jx|)^{-N}.
 \end{align}
\end{lemma}

\begin{proof}
For the first estimates \eqref{pointwiseA},  it seems to be
well-known. For instance, one can see \cite[Proposition
11.1]{LEM02}, \cite[Remark 2.2]{DDY} and \cite{LZ10,MYZ08} for the
details of the proof.

In order to prove the \eqref{pointwiseB},  recall the definitions of
$K^\alpha_{s,j}(t,x)$,  by scaling  we obtain
  that
  \begin{align}\label{DL-3.13A}
   \; \Big| \partial^\beta_x K^\alpha_{s,j}(t,x)\Big |
    &=\Big|\! \int_{\frac{3}{4}\le|\xi|\le \frac{8}{3}}e^{i\xi\cdot
    2^jx}e^{-t2^{2\alpha j}|\xi|^\alpha}2^{(s+n+|\beta|)j}\ \xi^\beta |\xi|^s\varphi(\xi)\ d\xi\Big|\nonumber \\
    & \!\le \!C\,2^{j(s+n+|\beta|)}\ e^{-ct2^{2\alpha j}}
  \end{align}
  and for any $N\ge1$,
 \begin{align}\label{DL-3.13B}
     &|2^j x|^{N} \Big|\partial^\beta_x K^\alpha_{s,j}(t,x)\Big |\nonumber\\
     &\le C_N \ 2^{j(s+n+|\beta|)}\sum_{|\mu|=N} \Big| \int_{\mathbb{R}^n}e^{i\xi\cdot
    2^jx} \partial^\mu_\xi \big(e^{-t2^{2\alpha j}|\xi|^{2\alpha}} \xi^\beta |\xi|^s\varphi(\xi)\big)\ d\xi \Big| \nonumber\\
    &\le{C'_N}\ 2^{j(s+n+|\beta|)}\int_{\frac{3}{4}\le|\xi|\le \frac{8}{3}}e^{-t2^{2\alpha j}|\xi|^{2\alpha}}\big(t^{N}2^{2\alpha Nj}+1\big)\big(\sum_{|\mu|\le N}|\partial^\mu_\xi\varphi(\xi)|\big) d\xi
\nonumber\\
     &\le{C''_N} 2^{j(s+n+|\beta|)} e^{-ct2^{2\alpha j}},
  \end{align}
 where we use the bound
  $e^{- t2^{2\alpha j}|\xi|^{2\alpha}}(t^{N}2^{2\alpha Nj}+1)
 \le C_N e^{-ct2^{2\alpha j}}$ as $\frac{3}{4}\le|\xi|\le \frac{8}{3}$.
 Combining \eqref{DL-3.13A} and
 \eqref{DL-3.13B}, we prove \eqref{pointwiseB}.
\end{proof}

   The following {\it endpoint bilinear estimate} follows by using
   Lemma \ref{Lem:2.2}
   and the classical Hardy-Littlewood-Sobolev (H-L-S) inequality (see
   \cite[p.\,353]{Stein:1993}).
 \begin{lemma}\label{Lem:2.3}
 Let $B(v,w)$ be defined as in \eqref{DLY-3}, $1<\alpha<\frac{5}{4}$ and $q_\alpha=\frac{3}{\alpha-1}$.
 Then for any $T>0$, there exists a positive constant $C_\alpha$ depending only on $\alpha$
 such that
 \begin{align}\label{DLY-2.3}
  &\|B(v,w)\|_{L^{q_\alpha}_xL^2_T\cap L^\infty_T\dot{F}^{-\alpha,2}_{q_\alpha}}\le {C_\alpha}
  \|v\|_{L^{q_\alpha}_xL^2_T}\|w\|_{L^{q_\alpha}_xL^2_T}.
 \end{align}
 \end{lemma}
 \begin{proof}
  To prove \eqref{DLY-2.3}, by applying Lemma \ref{Lem:2.2},  Young's inequality with respect to time
  variable,
  and H-L-S inequality with $1<\frac{q_\alpha}{2}<\infty$, $\frac{-1+\alpha}{3}=\frac{2}{q_\alpha}-\frac{1}{q_\alpha}$ to
  $B(v,w)$, we have
 \begin{align*}
    \|B(v,w)\|_{L^{q_\alpha}_xL^2_T}
  &=\Big\|\int_{0}^t e^{-(t-\tau)(-\Delta)^\alpha}\mathbb{P}\nabla\cdot(v\otimes{w})d\tau\Big\|_{{L}^{q_\alpha}_x{L}^{2}_T}\\
   &\lesssim
  \Big\|\!\!\int_{\mathbb{R}^3}\!\int_{0}^t\!\Big((t-\!\tau)^{\frac{1}{2\alpha}}\!+\!|x-\!y|\Big)^{-4}|(v\otimes{w})(y,\tau)|d\tau{dy}\Big\|_{{L}^{q_\alpha}_x{L}^{2}_T}\\
   &\lesssim
   \Big\|\!\!\int_{\mathbb{R}^3} |x-\!y|^{{-4+\alpha}}\|(v\otimes{w})(y,\cdot)\|_{L^1_T}{dy}\Big\|_{{L}^{q_\alpha}_x}\\
   &\lesssim
   \|v\otimes{w}\|_{{L}^{\frac{q_\alpha}{2}}_xL^1_T}.
 \end{align*}
 Next,  using  the boundedness of $\mathbb{P}$ in homogeneous
 space
 $\dot{F}^{-\alpha,2}_{q_\alpha}$, $\dot{F}^{0,2}_{q_\alpha}=L^{q_\alpha}_x$ 
 and Minkowski inequality,
 we get
 \begin{align}
    \|B(v,w)\|_{L^\infty_T\dot{F}^{-\alpha,2}_{q_\alpha}}& \lesssim
    \Big\|\int_{0}^t
     e^{-(t-\tau)(-\Delta)^\alpha}(v\otimes{w})d\tau\Big\|_{L^\infty_T\dot{F}^{1-\alpha,2}_{q_\alpha}}\nonumber\\
    &\lesssim\Big\|\int_{0}^t
     \Big(K_{1-\alpha}^\alpha(t-\tau,\cdot)\ast((v\otimes{w})(\tau,\cdot)\Big)(x) d\tau\Big\|_{L^\infty_TL^{{q_\alpha}}_x}\nonumber\\
    &\lesssim\Big\|\int_{0}^t
     \Big(K_{1-\alpha}^\alpha(t-\tau,\cdot)\ast((v\otimes{w})(\tau,\cdot)\Big)(x) d\tau\Big\|_{L^{q_\alpha}_xL^\infty_T},\label{DLY-2.4}
 \end{align}
 where $K^\alpha_{1-\alpha}(t,x)$ is the kernel of
 $e^{t\Delta}(-\Delta)^{\frac{1-\alpha}{2}}$.
  Making use of Lemma \ref{Lem:2.2}, we have
   \begin{align}\label{DLY-2.5}
     \sup_{t>0}|K_{1-\alpha}^\alpha(t, x\!-\!y)| \lesssim {|x-y|^{-4+\alpha}}.
   \end{align}
 Plugging \eqref{DLY-2.5} into \eqref{DLY-2.4}, using the H-L-S inequality and H\"older inequality, we get
   \begin{align*}
    (\ref{DLY-2.4})
     &\lesssim\Big\|\int_{\mathbb{R}^3}\frac{1}{|x-y|^{3-\alpha+1}}\|v\otimes{w}(y,\cdot)\|_{L^1_T}dy\Big\|_{L^{q_\alpha}_x}\\
     &\lesssim\Big\|\|v\otimes{w}\|_{L^1_T}\Big\|_{L^{\frac{q_\alpha}{2}}_x}\lesssim\|v\|_{L^{q_\alpha}_xL^2_T}\|w\|_{L^{q_\alpha}_xL^2_T}.
   \end{align*}
 Hence we finish the proof of \eqref{DLY-2.3}.
   \end{proof}
 The next lemma is about the equivalent definition and characterization
 of the Triebel-Lizorkin space which proof is given in the Appendix.
 \begin{lemma}\label{Lem:2.4}
 For any $1<\alpha<\frac{5}{4}$ and
 $q_\alpha=\frac{3}{\alpha-1}$, we have the following equivalent
 definition of $\dot{F}^{-\alpha,2}_{q_\alpha}$, i.e. any $f\in
 \dot{F}^{-\alpha,2}_{q_\alpha}$ satisfy
 \begin{align}\label{DLY:2.6}
   \Big\|\,f\Big\|_{\dot{F}^{-\alpha,2}_{q_\alpha}}\approx
   \Big(\int_{\mathbb{R}^3}(\int_{0}^\infty|e^{-t(-\Delta)^{\alpha}}f|^2dt)^{\frac{q_\alpha}{2}}dx\Big)^{{1}/{q_\alpha}}.
 \end{align}
 \end{lemma}
 \smallskip
  \begin{remark}\label{Rem:2.5}{\rm\hskip18cm

  (i) From \eqref{DLY:2.6}, we observe that given $u_0\in \dot{F}^{-\alpha,2}_{q_\alpha}$ and for any
  $\varepsilon>0$, there exists positive $T=T(u_0,\varepsilon)$ depending on the profile of $u_0$ such that
 $$\|e^{-t(-\Delta)^{\alpha}}u_0\|_{L^{q_\alpha}_xL^2_{T}}\le \varepsilon.$$
 One way to prove it is using splitting method or approximation by good functions.

 (ii) When $\alpha=1$, $\dot{F}^{-1,2}_\infty=BMO^{-1}$ and it has
 the following equivalent Carleson measure characterization which is closely related to \eqref{DLY:2.6} (see \cite{KOCT01}):
\begin{align*}
   \Big\|\,f\Big\|_{\dot{F}^{-1,2}_{\infty}}\approx
   \sup_{x\in\mathbb{R}^3,R>0}\Big(\frac{1}{B_R(x)}\int_{0}^{R^2}\int_{B_R(x)}|e^{t\Delta}f|^2dtdx\Big)^{\frac{1}{2}}.
 \end{align*}

  }
 \end{remark}

\subsection{The proof of Theorem 1.2}
 In order to prove well-posedness of the 3D gNS, we need to use the
 following Picard contraction principle, see for instance, \cite[Chapter 3.1, Lemma 4]{CAN04} and \cite{LEM02}.
 \vspace*{.4ex}
\begin{lemma}\label{Lem:2.6} Let
  $({X}, \|\cdot\!\|_{X})$  be an abstract Banach space and
  $B:{X}\times{X}\rightarrow{X}$ be a bilinear operator. If for any 
  $(u,v)\in{X}\times{X}$, there exists $c>0$ such that
 \begin{align*}
  \|B(u,v)\|_{X}\le{c}\|u\|_{X}\|v\|_{X},
 \end{align*}
  then for any $u_0$ satisfying
  $\|e^{t\Delta}u_{0}\|_{X}\!<\!{1}/{4c}$, the following system
  $u=e^{t\Delta}u_{0}+ B(u,u)$
  has a solution $u$ in $ {X}$. In particular, the solution is such that
  $\|u\|_{ {X}}\le{2}\|e^{t\Delta}u_{0}\|_{X}$
 and is the only one such that $\|u\|_{ {X}}<{1}/{2c}.$
\end{lemma}
  Now we are ready to prove the local and global well-posedness.

  \vspace*{1ex} \noindent\textbf{Proof of local
well-posedness}:\; Using Lemmas \ref{Lem:2.3}, \ref{Lem:2.4} and
\ref{Lem:2.6}, we prove that there exists a unique solution $u$ in a
closed ball in
 $L^{q_\alpha}_xL^2_T$
 since $\|e^{-t(-\Delta)^\alpha}u_0\|_{L^{\!q_\alpha}_x\!L^2_T}<\frac{1}{4c}$ if $u_0\in{\dot{F}^{-\alpha,2}_{q_\alpha}}$ and $T$ is small enough.
 Next we prove additional property of $u$ via Lemma \ref{Lem:2.3}, i.e.
 $u\!\in\!
 L^\infty_T\dot{F}^{-\alpha,2}_{q_\alpha}$.
 At last, following the standard dense argument we show that $u\in
 C([0,T];\dot{F}^{-\alpha,2}_{q_\alpha})$, see
 \cite{Kato:1984471} for details of Picard iteration arguments.

\vspace*{1ex} \noindent\textbf{Proof of global well-posedness}:
Noticing that bilinear estimates in Lemma \ref{Lem:2.3} can be
extended to $T=\infty$. However, in this case,
$\|e^{-t(-\Delta)^\alpha}u_0\|_{L^{q_\alpha}_xL^2_t}$ is not
necessarily small. Hence smallness condition is needed for global
well-posedness. Proof of global well-posedness follows in the
similar way.
 \subsection{Applications of other bilinear estimates}
  In this section, we consider many other bilinear estimates and give several applications to the analysis of well-posedness for the 3D gNS. Moreover, it is worth mentioning that these bilinear estimates are also very
  important to the ill-posedness of the 3D gNS.

  At first, we recall the endpoint bilinear estimates
  proved in Lemma \ref{Lem:2.3}:
  \begin{align}\label{A--0}
    \Big\|B(u,v)\Big\|_{L^{q_\alpha}_xL^2_T}
    &\!\!=\!\Big\|\!\int_0^t\!e^{-(t-\tau)(-\Delta)^\alpha}\mathbb{P}\nabla\cdot(u\otimes{v})d\tau\Big\|_{L^{\!q_\alpha}_xL^2_T}\!
    \lesssim \Big\||u||v|\Big\|_{{L^{\frac{q_\alpha}{2}}_xL^1_T}}\!.
  \end{align}

  It is clear that
  $L^{\frac{3}{2\alpha-1}}_xL^\infty_T$ is another endpoint space-time
  space for the 3D gNS equations. This kind of space-times space was first introduced by Calder\'{o}n
  \cite{CAL90} to study the incompressible
  Navier-Stokes equations. Applying maximal function theory and H-L-S inequality, we get
    $\|\displaystyle{\sup_{t>0}}|e^{-t(-\Delta)^\alpha}u_0|\|_{L^{\!\frac{3}{2\alpha-\!1}}_x}\lesssim\|u_0\|_{L^{\!\frac{3}{2\alpha-1}}_x}$ and
  \begin{align}\label{A--2}
  \Big\|B(u,v)\Big\|_{L^{\frac{3}{2\alpha-1}}_xL^\infty_T}
  &\lesssim
    \Big\|(-\!\Delta)^{\frac{1-2\alpha}{3}}\!(\||u||v|\|_{L^\infty_T})\Big\|_{L^{\frac{3}{2\alpha-1}}_x}
     \lesssim\Big\||u||v|\Big\|_{L^{\frac{3}{4\alpha-2}}_xL^\infty_T}.
   \end{align}
 Hence {\it{local} and global well-posedness of gNS for small} $L^{\!\frac{3}{2\alpha-1}}(\mathbb{R}^3)$ data follows from Picard contraction argument.
 Additionally, {\it local well-posedness of the 3D gNS for large
 $L^{\!\frac{3}{2\alpha-1}}(\mathbb{R}^3)$} data follows in the similar
 way as in \cite{Kato:1984471}.

\medskip
 For any $\frac{1}{q_{\theta}}=\frac{\theta}{q_\alpha}+\frac{(1-\theta)(2\alpha-1)}{3}$,
 $\frac{1}{p_\theta}=\frac{\theta}{2}+\frac{1-\theta}{\infty}$ and
 $0<\theta<1$, by Interpolating \eqref{A--0}--\eqref{A--2}, we have
 \begin{align}\label{A--5}
  \Big\|B(u,v)\Big\|_{L^{q_\theta}_xL^{p_\theta}_T}\lesssim
  \Big\||u||v|\Big\|_{L^{\frac{q_\theta}{2}}_xL^{\frac{p_\theta}{2}}_T}.
 \end{align}
 Thus \eqref{A--5} and (A.6) yield {\it well-posedness of the gNS
 for $\dot{F}^{-\theta\alpha,p_\theta}_{q_\theta}(\mathbb{R}^3)$} data.

 \medskip
 Notice that by applying the similar arguments as introduced by
 Fujita-Kato (cf. \cite{FUJK64}), one can prove local and global well-posedness for any
 $\dot{H}^{\frac{5-4\alpha}{2}}(\mathbb{R}^3)$ initial
 data. Thus interpolation between $\dot{F}^{-\alpha,2}_{q_\alpha}(\mathbb{R}^3)$ and
  $\dot{F}^{\frac{5-4\alpha}{2},2}_2(\mathbb{R}^3)$ yields {\it well-posedness} of
  the gNS in any Triebel-Lizorkin space $\dot{F}^{s_\theta,2}_{q_\theta}(\mathbb{R}^3)$ with
  $s_\theta=\theta(-\alpha)+\frac{(1-\theta)(5-4\alpha)}{2}$,
  $\frac{1}{q_\theta}=\frac{\theta}{q_\alpha}+\frac{1-\theta}{2}$
  and $0<\theta<1$ {\rm (\,see \cite[p.\,44]{Triebel:1978}\,).

 \medskip Recall the ideas that
 Cannone and Planchon used to prove the well-posedness of
 the N-S equations for data in the Besov
 space $B^{s,r}_{q}(\mathbb{R}^n)$. By making use of Lemma \ref{Lem:2.3}, it is easy to prove that
 for any $\frac{3}{2\alpha-1}<q<\frac{3}{\alpha-1}$ and
 $1<\alpha<\frac{5}{4}$,
 \begin{align}
  \sup_{t>0}t^{\frac{2\alpha-1-\frac{3}{q}}{2\alpha}}\|B(v,w)\|_{L^q_x}.
   &\lesssim \sup_{\tau>0}\,\Big(\tau^{\frac{2\alpha-1-\frac{3}{q}}{2\alpha}}\|v\|_{L^q_x}\Big)
   \sup_{\tau>0}\,\Big(\tau^{\frac{2\alpha-1-\frac{3}{q}}{2\alpha}}\|w\|_{L^q_x}\Big)\label{B-0}
 \end{align}
 Moreover, from \cite[Proposition 2.1]{MYZ08}, we have
 \begin{align}
 \sup_{t>0}t^{\frac{2\alpha-1-\frac{3}{q}}{2\alpha}}\|e^{-t(-\Delta)^\alpha}u_0\|_{L^q_x}\sim
 \|u_0\|_{\dot{B}^{1-2\alpha+\frac{3}{q},\infty}_{q}}.\label{B-1}
 \end{align}
 Combining \eqref{B-0}--\eqref{B-1}, we can prove well-posedness of
 the gNS in $\dot{B}^{1-2\alpha+\frac{3}{q},\infty}_q(\mathbb{R}^3)$
  with $\frac{1}{2\alpha-1}<q<q_\alpha=\frac{3}{\alpha-1}.$ Moreover,
 we have
 \begin{align}
 \label{B-3}
 \dot{B}^{1-2\alpha+\frac{3}{2},r\ge1}_2(\mathbb{R}^3)\hookrightarrow
 \dot{B}^{1-2\alpha+\frac{3}{q},\infty}_q(\mathbb{R}^3).
 \end{align}

 \medskip
   Next,  we shall consider the bicontinuity of $B(u,v)$ in some spaces which are not scale invariant. For instance, we prove
    the following bilinear estimate:
 \begin{align}
 \|B(u,v)&\|_{L^\infty_xL^2_T}
   \lesssim
   \|\!\int_{\mathbb{R}^3}\int_0^t\frac{|u(\tau,y)||v(\tau,y)|}{((t\!-\!\tau)^{\frac{1}{2\alpha}}\!+\!|x\!-\!y|)^{n+1}}d\tau dy\|_{L^\infty_xL^2_T}\nonumber\\
%
   &\lesssim
   T^{\frac{1}{2}}\Big\|\!\int_{|x-y|\ge1}\!\frac{\||u||v|\|_{L^1_T}}{|x\!-\!y|^{n+1}}dy\Big\|_{L^\infty_x}\!\!+
   \!\Big\|\!\int_{|x-y|\le1}\frac{\||u||v|\|_{L^1_T}}{|x\!-\!y|^{n-\alpha+1}}dy\Big\|_{L^\infty_x}\nonumber\\
   &\lesssim\max\{T^{\frac{1}{2}},\,1\}\Big\||u||v|\Big\|_{L^\infty_xL^1_T}.\label{A--4}
 \end{align}
 Interpolating \eqref{A--0} and \eqref{A--4}, then using (A.6), for any $q_\alpha<q<\infty$ and $T\lesssim1$, we
 have
 \begin{align}\label{D--3}
  \|B(u,v)\|_{L^q_xL^2_T}\lesssim
  \||u||v|\|_{L^{\frac{q}{2}}_xL^1_T},\quad
  \|e^{-t(-\Delta)^\alpha}u_0\|_{L^q_xL^2_t}\sim
  \|u_0\|_{\dot{F}^{-\alpha,2}_{q}},
 \end{align}
 which yields {\it local well-posedness of the 3D gNS with
 $\dot{F}^{-\alpha,2}_{q}$-valued initial data.} 

\section{ Analysis of ill-posedness }
 In this section, we will prove ``norm inflation" of the gNS
 in $\dot{F}^{-\alpha,r}_{q_\alpha}$ with
 $r>2$ and $q_\alpha=\frac{3}{\alpha-1}$.
 Following the
 ideas in \cite{BP08}, we rewrite
 the solution to the gNS equations as a  summation of the first approximation terms, the second
 approximation terms and remainder terms, i.e.
\begin{align}\label{eq:4.1}
 &u=u_1 \!- u_{2} + y,
 \end{align}
 where
 $u_1:=e^{-t(-\Delta)^\alpha}u_0:=S_tu_0$ and $u_{2}\!=B(u_1,u_1).$
 Moreover, the remainder terms
 satisfy the following integral equations:
 \begin{align}\label{eq:4.2}
 y &= G_0+ G_1 - G_2 ,
 \end{align}
       on $(0,\infty)$ with the
initial conditions $y(0)=0$,
\begin{align}\label{eq:4.3}
 \;\left\{\begin{aligned}
 G_2 & = B(y,y), \\
  G_1 & = B(y,u_2)+B(u_2,y)-B(y,u_1)-B(u_{1},y), \\
  G_0 & = B(u_{2}, u_1)+B(u_1,u_{2})- B(u_{2},u_{2}).
      \end{aligned}\right.
      \end{align}

 In the rest part of this section, we will establish the {\it
 a-priori} estimates for $u_0$, $u_{1}$, $u_2$ and $y$. Precisely,
 in Subsection 3.1, we construct some special initial data $u_0$; In Subsection 3.2 we estimate the small upper bounds of $u_0$ and $u_1$;
 In Subsection 3.3, we prove both upper bound and lower bound of $u_2$;
 In Subsection 3.4, we prove the upper bound of $y$; In
 Subsection 3.5, we complete the proof of Theorem \ref{Thm:1.3}.

 \subsection{Construction of initial data for the 3D gNS
 equations}
  For any fixed small number $\delta>0$, we define the initial data as follows:
 \begin{align}\label{eq:4.4}
   u_0
   &=\frac{Q}{\sqrt{\rho}}\sum_{s=1}^{\rho}
                              |k_s|^{\alpha}
                              \left\{\left(\!\!\!{\begin{array}{c}
                                  0        \\
                                  \frac{-\partial_3}{|k_s|}\psi_{k_s}         \\
                                   \frac{\partial_2}{|k_s|}\psi_{k_s}  \\
                                                 \end{array}}\!\!\!\right)
                       \!+\! \left(\!\!\!{ \begin{array}{c}
                                  0        \\
                                  \frac{-\partial_3}{|k_s|}\bar{\psi}_{k_s}           \\
                                  \frac{\partial_2}{|k_s|}\bar{\psi}_{k_s}  \\
                                                 \end{array}}\!\!\!\right)
                       \!+\!\left(\!\!\!\!{ \begin{array}{c}
                                  \frac{\partial_2}{|k_s|}{\psi_{k_s'}}   \\
                                  \frac{-\partial_1}{|k_s|}{\psi_{k_s'}}  \\
                                  0   \\
                                                 \end{array}}\!\!\!\right)
                       \!+\!\left(\!\!\!\!{ \begin{array}{c}
                                  \frac{\partial_2}{|k_s|}\bar{\psi}_{k_s'}   \\
                                  \frac{-\partial_1}{|k_s|}\bar{\psi}_{k_s'}   \\
                                  0   \\
                                                 \end{array}}\!\!\!\right)\right\}
                                               \nonumber\\
  :\!&\!=
  \frac{Q}{\sqrt{\rho}}\sum_{s=1}^\rho|k_s|^\alpha(\Psi_{k_s}\!
  +\!\bar{\Psi}_{k_s}\!+\!{\Phi_{k_s'}}\!+\!\bar{\Phi}_{k_s'})
 \end{align}
  with $\bar{\Psi}_{k_s}$, $\bar{\Phi}_{k_s'}$, $\bar{\psi}_{k_s}$ and $\bar{\psi}_{k_s'}$ being conjugate
  functions of $\Psi_{k_s}$,
  ${\Phi_{k_s'}}$, $\psi_{k_s}$ and ${\psi_{k_s'}}$, respectively.
 In addition, the parameters and auxiliary functions satisfy:
\begin{itemize}
 \item[({\bf H1})]
   $Q$, $\rho$ and  $m_0$ {will be chosen sufficiently large} {according to the size of} $\delta$ and
  \begin{align*}
    k_0=(0,2^{m_0},0),\;\;\; k_s=(0,2^{\frac{(s+1)(s+2m_0)}{2}},0),\;\;\; k_s'=(7,-2^{\frac{(s+1)(s+2m_0)}{2}},0),
  \end{align*}
   with $s=1,2,\cdots$ and $\rho$ will be specified in Lemma \ref{Lem:3.11}.
\item [({\bf H2})]
   $\psi(x)$, $\psi_{k_s}(x)$, ${\psi_{k_s'}}(x)$, $\bar{\psi}_{k_s}(x)$ and $\bar{\psi}_{k_s'}(x)$ satisfy:\;
 \begin{align*}
   & \widehat{\psi}(\xi)\!=\widehat{\psi}(|\xi|)\ge0,
     \;\;\;\;\;   \textrm{supp}\,\widehat{\psi}\subset
     {B}_{\frac{1}{4}}(0):=\{|\xi|<\frac{1}{4}\},\;\hskip2.2cm\\
   &
     \|\widehat{\psi}\|_{L^1_\xi}=1,\,\,\,\; \|\widehat{\psi}\|_{L^{\infty}_\xi}\!\sim \!\|\psi\|_{L^{q_\alpha}_x}\!\sim\!1, \;\;\;\widehat{\psi}_{k_s}(\xi)=\widehat{\psi}(\xi-k_s),
      \\
   &
   \widehat{\psi}_{k_s'}(\xi)=\widehat{\psi}(\xi-k_s'),\;\;
   \widehat{\bar{\psi}}_{k_s}(\xi)=\widehat{\psi}_{-k_s}(\xi),\;\;
     \widehat{\bar{\psi}}_{k_s'}(\xi)=\widehat{\psi}_{-k_s'}(\xi).
 \end{align*}
   \end{itemize}
 \smallskip
 \begin{remark} {\rm
    From hypothesis (H1)--(H2), we have the following observations:
  \begin{itemize}
  \vspace*{-1ex}
  \item[i)]
    $\psi(x)$ is real-valued and smooth,
    $\psi_{k_s}(x)=e^{ik_s\cdot x}\psi(x)$,
    $\bar{\psi}_{k_s}(x)=e^{i(-k_s)\cdot
    x}\psi(x)=\psi_{-k_s}(x).$
    As a consequence, $u_0(x)$ is real-valued, smooth and divergence free.
 \item[ii)]
    For any $(k,q)\!\in\!\mathbb{Z}^3\times[1,\infty]$, we have $\|\widehat{\psi}_{k}\|_{L^q_\xi}=\|\widehat{\psi}\|_{L^q_\xi}$ since $L^q_\xi$ is a
    shift invariant space.
    Making use of Hausdorff--Young's inequality and (H2), we get
 \begin{align}\label{eq:4.5}
   \|\Psi_{k_s}\|_{L^{q_\alpha}_x}+\|\bar{\Psi}_{k_s}\|_{L^{q_\alpha}_x}
       +\|{\Phi_{k_s'}}\|_{L^{q_\alpha}_x}+\|\bar{\Phi}_{k_s'}\|_{L^{q_\alpha}_x}
       \lesssim\|\widehat{\psi}\|_{L^{\!\frac{q_\alpha}{q_\alpha-1}}_\xi}\lesssim1.
  \end{align}
\item[iii)]
          \vspace*{-2ex}
 {\it Lacunarity of the sequence $\{|k_s|\}_{s=1}^\rho$}.
    According to the choices of $k_s$, $k_s'$ and $m_0$ in (H1), $\log_2|k_s|\!={\frac{(s+1)(s+2m_0)}{2}}$  is integer and $|k_s|< |k_s'|< |k_s|+7$.
    For any positive integer $s\in[1,\rho\,]\cap\mathbb{N}$, if we
    denote
    $j_s=\log_2|k_s|-\!1$, then $j_{s+1}-j_s\ge m_0$. For sufficiently large $m_0$ ( $\ge5$),  we have $7<\frac{1}{4}\,2^{m_0}=\frac{1}{4}{|k_0|}< \frac{1}{4}{|k_s|}$,
 \begin{align}\label{DLY-4.6}
     \frac{3}{4}\,2^{j_s}<|k_s|<|k_s'|<\frac{8}{3}\,2^{j_s}\!,\;
    \frac{3}{4}\,2^{j_s+1}<|k_s|<|k_s'|<\frac{8}{3}\,2^{j_s+1}\!.
 \end{align}
    For suitably large $m_0$ and $j\le j_s-1$ or $j\ge j_s+2$, from \eqref{DLY-4.6} we have
  \begin{align}\label{DLY-4.7}
     \Big(B_{\frac{1}{4}}(k_s)\cup{B}_{\frac{1}{4}}(k_s')\Big)\cap\Big\{\xi\in\mathbb{R}^3;\;
    \frac{3}{4}\,2^{j}<|\xi|<\frac{8}{3}\,2^{j}\Big\}=\emptyset.
  \end{align}
    From \eqref{DLY-1.4}--\eqref{DLY-1.5} and \eqref{eq:4.4}--\eqref{DLY-4.7}, for any
    $\ell\in\!\{0,1\}$ and $j_i\in\! \{j_1,\cdots,j_\rho\}$ we get
 \begin{align}\label{DLY-4.8}
    \Delta_{j_i+\ell} \sum_{s=1}^\rho \Psi_{k_s}  = \Delta_{j_i+\ell}\Psi_{k_i},\;\;
    \Delta_{j_i+\ell} \sum_{s=1}^\rho {\Phi_{k_s'}} = \Delta_{j_i+\ell}{\Phi}_{k_i'}.
 \end{align}
  Moreover, for any $j\!\in\!\mathbb{Z}\backslash\{j_1,j_1+1,j_2,j_2+1,\cdots,j_\rho,j_\rho+1\}$, from \eqref{DLY-4.7}
  we have
 \begin{align}\label{DLY-4.9}
    \text{$\Delta_j\sum_{s=1}^\rho \Psi_{k_s} = \sum_{s=1}^\rho \Delta_{j}\Psi_{k_s}\equiv0,\;\;\;
    \Delta_j\sum_{s=1}^\rho {\Phi_{k_s'}}= \sum_{s=1}^\rho \Delta_{j}{\Phi_{k_s'}}\equiv0.$}
 \end{align}
\end{itemize}
}
\end{remark}
\subsection{Estimates for initial data and the first approximation terms}
  In this subsection, we will estimate $u_0$ and $u_1=e^{-t(-\Delta)^\alpha}u_0$.
 \begin{lemma}\label{Lem:3.2}\! For any initial $u_0$ defined in \eqref{eq:4.4} and any $r\!\ge 2$ and $1<\alpha<\frac{5}{4}$, we
 obtain that
  \begin{align}\label{DLY-4.10}
   &{\|u_0\|_{\dot{F}^{-\alpha,r}_{q_\alpha}}\!\lesssim\!
    Q\rho^{\frac{1}{r}-\frac{1}{2}}},\;\;
    {\|u_1\|_{\dot{F}^{-\alpha,r}_{q_\alpha}}\!\lesssim\!
    Q\rho^{\frac{1}{r}-\frac{1}{2}}e^{-ct|k_0|^{2\alpha}}},
\end{align}
for some absolute  constant $c>0$.
 \end{lemma}
\begin{proof}
 We first deal with the cases $2\le r<\infty$. In view of the construction of $u_0$, it suffices to bound ${Q}{{\rho}^{-\frac{1}{2}}}{\sum_{s=1}^\rho} |k_s|^\alpha\Psi_{k_s}$
  and ${Q}{{\rho}^{-\frac{1}{2}}}{\sum_{s=1}^\rho}|k_s|^\alpha e^{-t(-\Delta)^\alpha}\Psi_{k_s}$.
 By Definition \ref{Def:1.1}, \eqref{DLY-4.6}--\eqref{DLY-4.9} and $|k_s|\sim2^{j_s}$,
  we obtain
  \begin{align}
    \|\sum_{s=1}^\rho
    |k_s|^\alpha\Psi_{k_s}\|_{\dot{F}^{-\alpha,r}_{q_\alpha}}
   &\!=\Big\|\ \|\big\{2^{-\alpha j}\Delta_j(\sum_{s=1}^\rho
    |k_s|^\alpha\Psi_{k_s})\big\}\|_{l^r}\ \Big\|_{L^{q_\alpha}_x}\nonumber\\
   &\!= \Big\|\Big(\sum_{i=1}^\rho \sum_{\ell=0,1}2^{-j_i\alpha r}\, |k_{i}|^{\alpha r}\, |\Delta_{j_i+\ell}
    \Psi_{k_{i}}|^r\Big)^{\frac{1}{r}} \Big\|_{L^{q_\alpha}_x}\nonumber\\
    &\!
    \lesssim \Big\|\Big(\sum_{i=1}^\rho\sum_{\ell=0,1}  |\Delta_{j_i+\ell}
    \Psi_{k_{i}}|^r\Big)^{\frac{1}{r}} \Big\|_{L^{q_\alpha}_x}.\label{DLY-4.11}
  \end{align}
 Note that $|\psi_{k_i}(x)|=|\psi(x)|$ for any $1\le i\le \rho$,  then for any $\ell=0,1$ and $1\le i \le \rho$, we have the following point-wise estimates
 \begin{align}
 |\Delta_{j_i+\ell}\Psi_{k_i}|\le \sum_{|\beta|=1}|(|k_i|^{-1}\partial^\beta\Delta_{j_i+\ell})\psi_{k_i}|
 \lesssim M\psi, \label{DLY-4.11A}
 \end{align}
  where $Mf$ denotes the standard Hardy-Littlewood maximal function of $f$.  Hence
 it follows from \eqref{DLY-4.11A}, Hardy-Littlewood theorem \cite[Chapter\,1,
 p.\,13]{Stein:1993} and $\|M\psi\|_{L^{q_\alpha}_x}\lesssim \|\psi\|_{L^{q_\alpha}_x}\lesssim1$
 that
 \begin{align}
   \eqref{DLY-4.11}&\!\lesssim
    \Big\|\Big(\sum_{i=1}^\rho\sum_{\ell=0,1} |M\psi|^r\Big)^{\frac{1}{r}} \Big\|_{L^{q_\alpha}_x}\lesssim
    \rho^{\frac{1}{r}}\|M\psi\|_{L^{q_\alpha}_x}\lesssim \rho^{\frac{1}{r}},\label{DLY-4.11B}
  \end{align}
which immediately concludes the desired estimates
$\|u_0\|_{\dot{F}^{-\alpha,r}_{q_\alpha}}\lesssim
Q\rho^{\frac{1}{r}-\frac{1}{2}}$.

Next to estimate $\|u_1\|_{\dot{F}^{-\alpha,r}_{q_\alpha}}$,
similarly as in (\ref{DLY-4.11}), we can obtain that
\begin{align}
    \|\sum_{s=1}^\rho
    |k_s|^\alpha e^{-t(-\Delta)^\alpha}\Psi_{k_s}\|_{\dot{F}^{-\alpha,r}_{q_\alpha}}
    \lesssim \Big\|\Big(\sum_{i=1}^\rho \sum_{\ell=0,1}|\Delta_{j_i+\ell}
    e^{-t(-\Delta)^\alpha}\Psi_{k_{i}}|^r\Big)^{\frac{1}{r}} \Big\|_{L^{q_\alpha}_x}.\label{DLY-4.11C}
  \end{align}
Note  that for $\ell=0,1$ and $1\le i\le\rho $,
\begin{align}\label{DLY-4.11F}
 |\Delta_{j_i+\ell}e^{-t(-\Delta)^\alpha}\Psi_{k_{i}}|\le\sum_{|\beta|=1}|k_i|^{-1} |\partial^\beta
 K_{j_i+\ell}^{\alpha}(t,\cdot)*\psi_{k_i}|,\end{align} where
$K_{j_i+\ell}^{\alpha}(t,x)$ is the kernel of the operator
$e^{-t(-\Delta)^\alpha}\Delta_{j_i+\ell}$. By the \eqref{pointwiseB}
in Lemma \ref{Lem:2.2}, we have the estimates
\begin{align}\label{DLY-4.11G}
 |\partial^\beta K_{j_i+\ell}^{\alpha}(t,x)|\lesssim e^{-ct|k_i|^{2\alpha}}2^{|\beta|j_i}\  2^{3j_i}(1+2^{j_i}|x|)^{-4},
 \end{align} for some $c>0$.
 Thus for any $\ell=0,1$ and $1\le i\le\rho $, it follows
from the (\ref{DLY-4.11F}) and \eqref{DLY-4.11G} that
\begin{align}
|\Delta_{j_i+\ell}e^{-t(-\Delta)^\alpha}\Psi_{k_{i}}|
&\lesssim  e^{-ct|k_i|^{2\alpha}} M\psi \label{DLY-4.11D}
\end{align}
Hence in view of (\ref{DLY-4.11C})--(\ref{DLY-4.11D}) and $|k_0|\le
|k_i|$, it immediately follows from Hardy-Littlewood maximal theorem
that
 \begin{align}
    \|u_1\|_{\dot{F}^{-\alpha,r}_{q_\alpha}}&\!\lesssim
    \frac{Q}{\sqrt{\rho}}e^{-ct|k_0|^{2\alpha}} \Big\|\Big(\sum_{i=1}^\rho\sum_{\ell=0,1}|(M\psi)(x)|^r\Big)^{\frac{1}{r}} \Big\|_{L^{q_\alpha}_x}
    \!\nonumber\\&\!
    \lesssim Q\rho^{\frac{1}{r}-\frac{1}{2}}e^{-ct|k_0|^{2\alpha}}.\label{DLY-4.11E}
  \end{align}
Thus we complete the proof for the cases $2\le r<\infty$.  Finally,
as $r=\infty$, similar to \eqref{DLY-4.11} and \eqref{DLY-4.11C}, we
can obtain that
$$ \Big\|\sum_{s=1}^\rho
    |k_s|^\alpha\Psi_{k_s}\Big\|_{\dot{F}^{-\alpha,\infty}_{q_\alpha}}\le \
     \Big\|\sup_{1\le s\le \rho} (|\Delta_{j_s}\Psi_{k_s}|+|\Delta_{j_s+1}\Psi_{k_s}|)\Big\|_{L^{q_\alpha}_x}$$
and
$$ \Big\|\sum_{s=1}^\rho
    |k_s|^\alpha e^{-t(-\Delta)^\alpha}\Psi_{k_s}\Big\|_{\dot{F}^{-\alpha,\infty}_{q_\alpha}}\le \
     \Big\|\sup_{1\le s\le \rho} \sum_{\ell=0,1}|\Delta_{j_s+\ell}\,e^{-t(-\Delta)^\alpha}\Psi_{k_s}|\Big\|_{L^{q_\alpha}_x}$$
Hence by the estimates \eqref{DLY-4.11A} and \eqref{DLY-4.11D} we
can immediately conclude the desired bounds of $u_0$ and $u_1$ for
the case $r=\infty$.
\end{proof}
\vskip0.3cm
\begin{remark}\label{Rem:3.3} {\rm From \eqref{DLY-4.10}, for any given $\delta>0$, if $r>2$, then there exists sufficiently large
 $\rho$ such that $\|u_0\|_{\dot{F}^{-\alpha,r>2}_{q_\alpha}}\le \delta$ and $\|u_1\|_{\dot{F}^{-\alpha,r>2}_{q_\alpha}}\le \delta$.
  Similarly, {\it it is also easy to prove that for any $q\in(q_\alpha,
  \infty]$,}
  \begin{align}
  \label{D-0}\|u_0\|_{\dot{F}^{-\alpha,r>2}_{q}}\lesssim Q\rho^{\frac{1}{r}-\frac{1}{2}},\quad
   \|u_1\|_{\dot{F}^{-\alpha,r>2}_{q}}\lesssim  Q\rho^{\frac{1}{r}-\frac{1}{2}}. \end{align}
  Therefore, for sufficiently large $\rho$,
  $\|u_0\|_{\dot{F}^{-\alpha,r>2}_{q}}\le\delta$
  and
  $\|u_1\|_{\dot{F}^{-\alpha,r>2}_{q}}\le\delta$.
   }
 \end{remark}

\medskip
 \begin{lemma}\label{Lem:3.4}
  For any $T>0$, $u_0$ and $u_1$ given in \eqref{eq:4.4} and
  \eqref{eq:4.1}, we obtain that
 \begin{align}\nonumber
   \|u_1\|_{L^{q_\alpha}_xL^2_T}\lesssim
   \frac{Q}{\sqrt{\rho}}(T^{\frac{1}{2}}|k_{\rho-N_0}|^\alpha+\sqrt{N_0}),
 \end{align}
for any $0\le N_0\le \rho$. In particular,
$\|u_1\|_{L^{q_\alpha}_xL^2_T}\rightarrow 0$ as
 $N_0=0$ and $T\rightarrow 0$.
 \end{lemma}
 \begin{proof}
 {\it From now on, we let $S_tu_0=e^{-t(-\Delta)^\alpha}u_0=u_1$.}
 By the construction of initial data $u_0$, it suffices to
 estimate $\displaystyle{\sum_{s=1}^\rho}|k_s|^\alpha S_t\Psi_{k_s}$.
 Similar to \eqref{DLY-4.11}), we
 have
\begin{align}
  \| \sum_{s=1}^\rho |k_s|^\alpha S_t\Psi_{k_s}\|_{L^{q_\alpha}_xL^2_T}
  &\lesssim\,\| \sum_{s=1}^\rho |k_s|^\alpha S_t\Psi_{k_s}\|_{L^2_T\dot{F}^{0,2}_{q_\alpha}}\nonumber\\
  &\lesssim  \Big\||M\psi|\Big(
     \sum_{s=1}^\rho\sum_{\ell=0}^1
     |k_s|^{2\alpha}e^{-ct|k_s|^{2\alpha}}\Big)^{\frac{1}{2}}\Big\|_{L^2_TL^{q_\alpha}_x}\nonumber\\
     &\lesssim
     \Big\|\Big(
     \sum_{s=1}^\rho
     |k_s|^{2\alpha}e^{-ct|k_s|^{2\alpha}}\Big)^{\frac{1}{2}}\Big\|_{L^2_T} \Big\|M\psi\Big\|_{L^{q_\alpha}_x}\nonumber\\
     &\lesssim
      \Big( \sum_{s=1}^\rho  \int_0^T\, |k_s|^{2\alpha}e^{-ct|k_s|^{2\alpha}}dt\Big)^{\frac{1}{2}}.\label{DLY-4.12}
  \end{align}
 Notice that
  $$ \int_0^T\, |k_s|^{2\alpha}e^{-ct|k_s|^{2\alpha}}dt \lesssim  \min\{1,
  \,T|k_s|^{2\alpha}\}
 \text{ and
 }\displaystyle{\sum_{s=1}^{\rho-N_0}}|k_s|^{2\alpha}\!\lesssim\!|k_{\rho-N_0}|^{2\alpha},$$
 then we get
 \begin{align}
  \eqref{DLY-4.12}&\lesssim
   \Big(\sum_{s=1}^{\rho-N_0} T|k_s|^{2\alpha}\!+\!\!\sum_{s=\rho-N_0+1}^\rho\!\!1\Big)^{\frac{1}{2}}
  \lesssim  ({T}^{\frac{1}{2}}|k_{\rho-N_0}|^{\alpha}\!+\!\sqrt{N_0}).\label{DLY-4.13A}
  \end{align}
  Similar to \eqref{DLY-4.13A}, we can obtain the desired estimate.
\end{proof}

 {\rm By checking the estimates \eqref{DLY-4.13A} for the case $N_0=\rho$
again,  we know that the best upper bound of
 $\|u_1\|_{L^{q_\alpha}_xL^2_T}$ is actually $c\,Q$, which is not good enough to bound the remainder
 $y(t,x)$ (below). 
  Recalling the idea of estimating remainder term by means of bi-continuity of bilinear
 operator $B(u,v)$ in Lemma \ref{Lem:2.3},
 it is also natural to hope that nonlinear terms of $y(t,x)$ are smaller.
 Therefore, we need to analyze how $y$ evolve in different time scales
 and see their contributions by using the
 time-step-division method introduced by
 Bourgain-Pavlovi\'c in \cite{BP08} to prove ill-posedness of the 3D incompressible
 Navier-Stokes equations.
 Let
 \begin{align}\label{DLY-4.14}
 |k_\rho|^{-2\alpha}=T_0 < T_1 < T_2 < \cdots < T_\beta=|k_0|^{-2\alpha},
 \end{align}
  where $\beta = Q^3$, $T_\sigma = |k_{\rho_\sigma}|^{-2\alpha}$, $\rho_\sigma = \rho -\sigma Q^{-3}\rho$ and
  $\sigma =0, 1, 2,\cdots,\beta$.

\medskip
\begin{lemma}\label{Lem:3.6}
 Assume that $u_0$ satisfy \eqref{eq:4.4} and
 $u_1=S_tu_0=e^{-t(-\Delta)^{\alpha}}u_0$, we have
 \begin{align}
  \label{DLY-4.15}
  \Big\|\,u_1\Big\|_{L^{q_\alpha}_xL^2_{[T_\sigma,T_{\sigma+1}]}}:=\Big\|\,u_1\chi_{_{[T_\sigma,T_{\sigma+1}]} }(t)\Big\|_{L^{q_\alpha}_xL^2_{T_{\sigma+\!1}}}
  \lesssim \frac{Q}{\sqrt{\rho}}(1+\sqrt{\rho Q^{-3}}).
 \end{align}
\end{lemma}
\begin{proof}
  Recall that $u_1=S_tu_0$ and
  $u_0=\frac{Q}{\sqrt{\rho}}\displaystyle{\sum_{s=1}^\rho|k_s|^\alpha}{S}_t(\Psi_{k_s}\!+\!\bar{\Psi}_{k_s}\!+\!{\Phi_{k_s'}}\!+\!\bar{\Phi}_{k_s'})$.
  It also suffices to estimate $\sum_{s=1}^\rho|k_s|^\alpha{S}_t\Psi_{k_s}$. Similar to \eqref{DLY-4.12}, by using
  $$ \int_{T_\alpha}^{T_{\alpha+1}}\, |k_s|^{2\alpha}e^{-ct|k_s|^{2\alpha}}dt \lesssim  \min\Big\{T_{\alpha+1}|k_s|^{2\alpha},\,\;1,\,\;e^{-cT_\alpha |k_s|^2}\Big\},$$
 and \eqref{DLY-4.14} we obtain that
$$ \text{$\sum_{s=1}^{\rho_{\alpha+1}} T_{\alpha+1}|k_s|^2\lesssim1$,
 $\sum_{s=\rho_{\alpha+1}+1}^{\rho_\alpha-1}1\lesssim \rho{Q}^{-3}$, $\sum_{s=\rho_\alpha}^\rho
 e^{-T_\alpha|k_s|^2}\lesssim1$}$$
 and
  \begin{align}\nonumber
 \|u_0\chi_{_{[T_\alpha,T_{\alpha+1}]}}(t)\|_{L^{q_\alpha}_xL^2_{T_{\alpha+1}}}
 \!\!
 &\lesssim  \frac{Q}{\sqrt{\rho}}\Big(\sum_{s=1}^\rho  \int_{T_\alpha}^{T_{\alpha+1}}\!\!\!e^{-ct|k_s|^{2\alpha}}|k_s|^{2\alpha}dt \Big)^{\frac{1}{2}}%
 \\
 &\lesssim \frac{Q}{\sqrt{\rho}}(1+\sqrt{\rho Q^{-3}}).
  \end{align}
 Thus we prove \eqref{DLY-4.15}.
 \end{proof}
 The following result is a consequence of $\sum_{s=1}^\rho e^{-T_\beta|k_s|^2}\lesssim1$ and Lemma \ref{Lem:3.6}.
\begin{corollary}\label{Cor:3.7}
 For any $T>T_{\beta}=|k_0|^{-2}=2^{-2m_0}$, we have
 \begin{align}
  \label{DLY-4.18}
  \Big\|\,u_1\Big\|_{{L^{q_\alpha}_xL^2_{{[T_\beta,T]}}}}=\Big\|u_1\chi_{_{[T_\beta,T]} }(t)\Big\|_{{L^{q_\alpha}_xL^2_{T}}}\lesssim
  \frac{Q}{\sqrt{\rho}}.
   \end{align}
\end{corollary}
\subsection{Estimates for the second approximation terms}
 We start this subsection by making some preliminary calculations.
 Recall that $u_1(\tau)=S_\tau
 u_0=e^{\tau\Delta}u_0$. In order to study the bilinear form  $u_2=B(u_{1},u_{1})$, from the construction of initial data $u_0$, we
 first split the {\it second approximation terms} $u_2$ into
 $$u_{2}={u_{2,0}}+u_{2,1}+u_{2,2},$$
 where
    \begin{align}
    u_{2,0}&=\frac{Q^2}{\rho}\sum_{s=1}^\rho \int_0^t  |k_s|^{2\alpha}S_{t-\tau}\mathbb{P}\,F_{s}\,d\tau,
             \label{DLY-4.19}\\
   u_{2,1}&=\frac{Q^2}{\rho}\sum_{s=1}^\rho \int_0^t  |k_s|^{2\alpha}S_{t-\tau}\mathbb{P}\,G_{s}\,d\tau,
             \label{DLY-4.20}\\
   u_{2,2}&=\frac{Q^2}{\rho}\sum_{s=1}^\rho\sum_{l\ne{s}} \int_{0}^t |k_s|^\alpha|k_l|^\alpha{S_{t-\tau}}\mathbb{P}\,{H}_{s,l}\,d\tau    \hskip2cm
   \label{DLY-4.21}
 \end{align}
 and
\begin{align}\label{DLY-4A}
 \left\{\begin{aligned}
 F_s\,& ={S_\tau}{\Phi}_{k_s'} \cdot \nabla{S_\tau}{\Psi}_{k_s}  +{S_\tau}\bar{\Phi}_{k_s'} \cdot \nabla{S_\tau}\bar{\Psi}_{k_s}
       +{S_\tau}{\Psi}_{k_s} \cdot \nabla{S_\tau}{\Phi}_{k_s'} \\
      &+{S_\tau}\bar{\Psi}_{k_s} \cdot \nabla{S_\tau}\bar{\Phi}_{k_s'} +{S_\tau}{\Phi}_{k_s'} \cdot \nabla{S_\tau}\bar{\Phi}_{k_s'} +{S_\tau}\bar{\Phi}_{k_s'} \cdot \nabla{S_\tau}{\Phi}_{k_s'}\\
      & +{S_\tau}{\Psi}_{k_s} \cdot \nabla{S_\tau}\bar{\Psi}_{k_s} +{S_\tau}\bar{\Psi}_{k_s} \cdot \nabla{S_\tau}{\Psi}_{k_s},\\
 G_s\,&= {S_\tau}{\Phi}_{k_s'} \cdot \nabla{S_\tau}\bar{\Psi}_{k_s} +  {S_\tau}\bar{\Phi}_{k_s'} \cdot \nabla{S_\tau}{\Psi}_{k_s}
       +{S_\tau}{\Psi}_{k_s} \cdot \nabla{S_\tau}\bar{\Phi}_{k_s'} \\
      &+ {S_\tau}\bar{\Psi}_{k_s} \cdot \nabla{S_\tau}{\Phi}_{k_s'} + {S_\tau}{\Phi}_{k_s'} \cdot \nabla{S_\tau}{\Phi}_{k_s'} + {S_\tau}\bar{\Phi}_{k_s'} \cdot \nabla{S_\tau}\bar{\Phi}_{k_s'}\\
      &+ {S_\tau}{\Psi}_{k_s} \cdot \nabla{S_\tau}{\Psi}_{k_s} +  {S_\tau}\bar{\Psi}_{k_s} \cdot \nabla{S_\tau}\bar{\Psi}_{k_s}), \\
 H_{s,l}&= S_\tau(\Psi_{k_s}\!+\bar{\Psi}_{k_s}\!+{\Phi_{k_s'}}\!+\bar{\Phi}_{k_s'}) \cdot
 \nabla{S_\tau}(\Psi_{k_l}\!+\bar{\Psi}_{k_l}\!+{\Phi}_{k_l'}\!+\bar{\Phi}_{k_l'}).
 \end{aligned}\right.
 \end{align}

 \vskip0.3cm
 \begin{remark}
 {\rm Noticing that $u_{2,\ell}$ $(\ell=0,1,2)$ are {\it
 real-valued smooth}
 functions since they are summations of conjugated smooth functions.
 \, From \eqref{eq:4.4} and (H1)--(H2), we observe
 that $\widehat{F}_s$, $\widehat{G}_s$ and $\widehat{H}_{s,l}$ are
 {\it purely imaginary} since $\widehat{S_\tau\nabla{\Phi}_{k_l'}}$ is real-valued
 and $\widehat{S_\tau\Psi_{k_s}}$
   is purely imaginary.}

  \medskip
  {\rm According to the different frequency interactions, we decompose the second approximation terms $u_{2}$ into three parts which are given
  in \eqref{DLY-4.19}--\eqref{DLY-4.21}.
  Precisely,
  \begin{itemize}
  \item $u_{2,1}$ represents the {\it high-high to high} frequency
  interactions. This is always the best one, see Lemma \ref{Lem:3.9}.
  \item $u_{2,2}$ represents the {\it high-low to high} and {\it low-high to
  high} frequency interactions. Usually, these two kinds of frequency
  interactions can be well-controlled, see Lemma \ref{Lem:3.10}
  \item $u_{2,0}$ represents the
   {\it high-high to low} frequency interactions. This kind of interaction
   is always the worst one. In this paper, we will explore
   this part and gain the desired lower bound. Proof of the lower bound is a little complicated, thus we would like to give the
   details in the end of this subsection, see Lemma \ref{Lem:3.11}.
  \end{itemize}

  }
  \end{remark}
 Now we prove the following estimates for $u_{2,1}$.
  \begin{lemma}\label{Lem:3.9}
 For any $1<\alpha<\frac{5}{4}$, $r\ge2$ and $q_\alpha=\frac{3}{\alpha-1}$, we
 have
\begin{align}\label{DLY-4.22}
&\|u_{2,1}\|_{L^\infty_T\dot{F}^{-\alpha,r}_{q_\alpha}}
+\|u_{2,1}\|_{L^{q_\alpha}_xL^2_T}\lesssim \frac{Q^2}{\rho}.
\end{align}
\end{lemma}
\begin{proof}

 First deal with the norm $\|u_{2,1}\|_{L^\infty_T\dot{F}^{-\alpha,r}_{q_\alpha}}$. Recall that from the initial value
 construction in
 \eqref{eq:4.4},
  $$\text{$\supp\widehat{\Psi}_{k_s}\!\!\subset\!\!{B}_{\frac{1}{4}}(k_s)$,\,
  $
  \supp\widehat{\Phi}_{k_s'}\!\!\subset\!\!{B}_{\frac{1}{4}}(k_s')$,\,
  $
  \supp\widehat{\bar{\Psi}}_{k_s}\!\!\subset\!\!{B}_{\frac{1}{4}}(-k_s)$,\,
  $ \supp\widehat{\bar{\Phi}}_{k_s'}\!\!\subset\!\!{B}_{\frac{1}{4}}(-k_s').
  $ }$$
Hence for any $1\le s\le \rho$, by \eqref{DLY-4.6}--\eqref{DLY-4.9},
we have 
 \begin{align}\label{DLY-4.23}
    G_{s}=(\Delta_{j_s+1}+\Delta_{j_s+2})G_{s},\ \  \Delta_{j}G_{s}=0,\ j\neq j_s+1 \text{ and } j\ne j_s+2.
 \end{align}
 which shows that $\{G_s\}_{s=1}^\rho$ is also a lacunary sequence.
 Therefore, by using boundedness of $\mathbb{P}$ and $|k_s|\sim 2^{j_s}$,  similar to \eqref{DLY-4.11},
 we obtain that
  \begin{align}\label{DLY-4.24}
   \|u_{2,1}\|_{L^\infty_T\dot{F}^{-\alpha,r}_{q_\alpha}}
    &\lesssim \frac{Q^2}{\rho}\ \Big\|\ \big\|\Big\{2^{-\alpha j}\Delta_j(\sum_{s=1}^{\rho}\!\int_0^t|k_s|^{2\alpha}S_{t-\tau}G_{s}
     d\tau)\Big\}\big\|_{l^r}\ \Big\|_{L^\infty_TL^{q_\alpha}_x}\nonumber\\
   &\!\lesssim \frac{Q^2}{\rho}\Big\|\Big(\sum_{s=1}^\rho \sum_{\ell=1,2}|k_{s}|^{\alpha r }\ \big|\int_0^tS_{t-\tau}\Delta_{j_s+\ell}G_{s}
   d\tau\big|^r\Big)^{\frac{1}{r}} \Big\|_{L^\infty_T
   L^{q_\alpha}_x}.
   \end{align}
 From Lemma \ref{Lem:2.2} and similar to \eqref{DLY-4.11D}, we get
  \begin{align}\label{DLY-4.25A}
  |(S_{t-\tau}\Delta_{j_s+\ell}G_{s})(x)|
  &\lesssim e^{-c(t-\tau)|k_s|^{2\alpha}} M_{|G_s|}(x).
  \end{align}
 where $M_{|G_s|}(x)$ denotes Hardy-Littewood maximal function of $|G_s|$.

Note that $\Psi_{k_s}=(\Delta_{j_s}+\Delta_{j_s+1})\Psi_{k_s}$ and
$\Psi_{k'_s}=(\Delta_{j_s}+\Delta_{j_s+1})\Psi_{k'_s}$. Hence by the
\eqref{DLY-4.11D}, we have ( if necessary, $c$ can be adjusted to be
smaller)
\begin{align}\label{DLY-4.25B}
|S_\tau\Phi_{k'_s}|+|S_\tau\Psi_{k_s}|+|S_\tau\bar{\Phi}_{k'_s}|+|S_\tau\bar{\Psi}_{k_s}|\lesssim
e^{-\frac{c}{2}\tau|k_s|^{2\alpha}}  M\psi
 \end{align}
and
 \begin{align} \label{DLY-4.25C}
 |\nabla{S_\tau}\Phi_{k'_s}|+|\nabla{S_\tau}\Psi_{k_s}|+ |\nabla{S_\tau}\bar{\Phi}_{k'_s}|+|\nabla{S_\tau}\bar{\Psi}_{k_s}|
\lesssim  |k_s|e^{-\frac{c}{2}\tau|k_s|^{2\alpha}} M \psi .
 \end{align}
Thus  \begin{align}
  |G_s|
  \lesssim
  |S_\tau{\Phi_{k_s'}}\cdot\nabla{S_\tau}\bar{\Psi}_{k_s}|+\cdots
 \lesssim |k_s|\ {e}^{-c\tau|k_s|^{2\alpha}}\ (M \psi)^2.\label{DLY-4.25}
 \end{align}
Let $\theta(x)=((M \psi)(x))^2$. Then it follows from
\eqref{DLY-4.25A} and \eqref{DLY-4.25} that
\begin{align}\label{DLY-4.25D}
  |S_{t-\tau}\Delta_{j_s+\ell}G_{s}(x)|
  \lesssim |k_s|  e^{-ct|k_s|^{2\alpha}}(M\theta)(x)
  \end{align}
 where $(M\theta)(x)$ is the maximal function of $\theta(x)$. Hence for any $1\le s\le\rho$ and $\ell=1,2$, we have
 \begin{align}\label{DLY-4.25E}
  \int_0^t|S_{t-\tau}\Delta_{j_s+\ell}G_{s}(x)|d\tau &\! \lesssim t|k_s|  e^{-ct|k_s|^{2\alpha}} (M\theta)(x) \lesssim |k_s|^{1-2\alpha} (M\theta)(x).
 \end{align}
 Now plugging \eqref{DLY-4.25E} into \eqref{DLY-4.24} and recalling that $\alpha>1$, then by Hardy-Littlewood maximal theorem we immediately obtain that
 \begin{align}\label{DLY-4.25F}
  \|u_{2,1}\|_{L^\infty_T\dot{F}^{-\alpha,2}_{q_\alpha}}&\!\lesssim
  \frac{Q^2}{\rho}\Big\|\Big(\sum_{s=1}^\rho |k_{s}|^{ r(1-\alpha)} \Big)^{\frac{1}{r}} M\theta\Big\|_{L^{q_\alpha}_x}
  \lesssim \frac{Q^2}{\rho} \|\theta\|_{L^{q_\alpha}_x}\lesssim \frac{Q^2}{\rho},
 \end{align}
 which concludes the desired bound of $\|u_{2,1}\|_{L^\infty_T\dot{F}^{-\alpha,2}_{q_\alpha}}$.

 For the norm of $\|u_{2,1}\|_{L^{q_\alpha}_xL^2_T}$. Similar to \eqref{DLY-4.24}, by using \eqref{DLY-4.25D} we have
 \begin{align}\label{DLY-4.26}
   \|u_{2,1}\|_{L^{q_\alpha}_xL^2_T}&\!\le  \|u_{2,1}\|_{L^2_T L^{q_\alpha}_x}=\|u_{2,1}\|_{L^2_T\dot{F}^{0,2}_{q_\alpha}}\nonumber\\
   &\! \lesssim \frac{Q^2}{\rho}\Big\|\Big(\sum_{s=1}^\rho \ \big| |k_{s}|^{2\alpha+1}t e^{-ct|k_s|^{2\alpha}} \big|^2\Big)^{\frac{1}{2}} (M\theta)(x)\Big\|_{L^2_T L^{q_\alpha}_x}\nonumber\\
    &\! \lesssim \frac{Q^2}{\rho}\ \Big(\sum_{s=1}^\rho \  \int_0^\infty\big||k_{s}|^{2\alpha+1}t e^{-ct|k_s|^{2\alpha}} \big|^2dt\Big)^{\frac{1}{2}} \ \|M\theta\|_{L^{q_\alpha}_x}\nonumber\\
     &\! \lesssim \frac{Q^2}{\rho}\ \Big(\sum_{s=1}^\rho \ k_s^{2(1-\alpha)}\Big)^{\frac{1}{2}} \lesssim \frac{Q^2}{\rho}.
 \end{align}
Thus combining \eqref{DLY-4.25F} and \eqref{DLY-4.26}, we finish the
proof of \eqref{DLY-4.22}.
\end{proof}

 Next we estimate $u_{2,2}$.
 \begin{lemma}\label{Lem:3.10}
 For any $1<\alpha<\frac{5}{4}$, $r\ge2$ and $q_\alpha=\frac{3}{\alpha-1}$, we
 obtain that
\begin{align}\label{DLY-4.27}
 &\|u_{2,2}\|_{L^\infty_T\dot{F}^{-\alpha,r}_{q_\alpha}}+\|u_{2,2}\|_{L^{q_\alpha}_xL^2_T}
  \lesssim \frac{Q^2}{{\rho}}.
\end{align}
\end{lemma}
 \begin{proof}
   In order to estimate $u_{2,2}$, we first rewrite it as follows:
  \begin{align}\label{DLY-4.28}
    u_{2,2}=&\frac{Q^2}{\rho}\sum_{s=1}^\rho |k_s|^\alpha(\int_0^t (\sum_{l<s}|k_l|^\alpha S_{t-\tau}\mathbb{P}H_{s,l})d\tau \nonumber\\
     &\hskip2cm + \frac{Q^2}{\rho}\sum_{l=1}^\rho|k_l|^\alpha{\int_0^t (\sum_{s<l}|k_s|^\alpha S_{t-\tau}\mathbb{P}H_{s,l})d\tau}.
  \end{align}
  Recall from the initial value construction \eqref{eq:4.4}, %
 we obtain that
  \begin{align*}
  \supp\widehat{H}_{s,l}\subset {B}_{\frac{1}{2}}(\pm k_s\!\pm\!k_l)
  \cup{B}_{\frac{1}{2}}(\pm k_s\pm k_l')
  \cup{B}_{\frac{1}{2}}(\pm k_s'\pm k_l)
  \cup{B}_{\frac{1}{2}}(\pm k_s'\pm k_l').
  \end{align*}
 As a result, we have
  \begin{align}\label{DLY-4.29A}
  \left\{\begin{aligned}
  &\sum_{{l<s}}H_{s,l}=\sum_{l\le s-1}
    (\Delta_{j_s}+\Delta_{j_s+1})H_{s,l}, \; \Delta_jH_{s,l}=0,\;
    j\ne j_s \text{ and } j\ne j_s+1,\\
   & \sum_{{s<l}}H_{s,l}=\sum_{s\le l-1}
    (\Delta_{j_l}+\Delta_{j_l+1})H_{s,l}, \; \Delta_jH_{s,l}=0,\;
    j\ne j_l \text{ and } j\ne j_l+1.
  \end{aligned}\right.
  \end{align}
 By checking the estimates \eqref{DLY-4.25B}--\eqref{DLY-4.25C} and
 recalling that $\theta(x)=(M\psi)^2(x),$
 for any $(s,l)\in (\mathbb{N}\cap[1,\rho])\times (\mathbb{N}\cap[1,\rho])$, we get
 \begin{align}\label{DLY-4.29C}
    |H_{s,l}|\lesssim |k_l|e^{-c\tau\max\{|k_s|^{2\alpha},\, |k_l|^{2\alpha}\}
    } (M\theta)(x)
 \end{align}
 and for $\kappa=0,1$, similar to \eqref{DLY-4.25A}, by using \eqref{DLY-4.29C} we get
 \begin{align}\label{DLY-4.29D}
  |S_{t-\tau}&\Delta_{\kappa+j_{\max\{s,\,l\}}}H_{s,l}|(x) \lesssim  |k_l|e^{-ct\max\{|k_s|^{2\alpha},\, |k_l|^{2\alpha}\}
    } (M\theta)(x).
 \end{align}
 In order to deduce the \eqref{DLY-4.27}, it suffices to estimate the first part of the \eqref{DLY-4.28}:
 $$\Big\|\frac{Q^2}{\rho}\sum_{s=1}^\rho\! \int_0^t\!|k_s|^\alpha(\sum_{l\le s-1}\!|k_l|^\alpha S_{t-\tau}\mathbb{P}H_{s,l})d\tau\Big\|_{L^\infty_T\dot{F}^{-\alpha,r}_{q_\alpha}\cap
 L^{q_\alpha}_xL^2_T}.$$
 Similar to \eqref{DLY-4.24}, \eqref{DLY-4.25F} and \eqref{DLY-4.29D}, by using boundedness of $\mathbb{P}$ in $\dot{F}^{-\alpha,r}_{q_\alpha}$, we have
 \begin{align}\label{DLY-4.29E}
 \|\frac{Q^2}{\rho}\sum_{s=1}^\rho\!
 \int_0^t& |k_s|^\alpha(\sum_{l\le s-1}\!|k_l|^\alpha
 S_{t-\tau}\mathbb{P}H_{s,l})d\tau\|_{L^\infty_T\dot{F}^{-\alpha,r}_{q_\alpha}}\nonumber\\
 &\lesssim\frac{Q^2}{\rho}\Big\|
  \Big(
   \sum_{s=1}^\rho\sum_{\kappa=0,1}
   (\sum_{l=1}^{s-1}|k_l|^{\alpha}|\int_0^tS_{t-\tau}\Delta_{j_s+\kappa}H_{s,l}d\tau|)^r
   \Big)^{\frac{1}{r}}\Big\|_{L^\infty_TL^{q_\alpha}_x}\nonumber\\
 &\lesssim\frac{Q^2}{\rho}\Big\|
  M\theta \Big(
   \sum_{s=1}^\rho
    (\sum_{l=1}^{s-1}|k_l|^{\alpha}|k_s|^{1-2\alpha})^r
   \Big)^{\frac{1}{r}} \Big\|_{L^{q_\alpha}_x}\nonumber\\
 &\lesssim\frac{Q^2}{\rho}
  \Big(
   \sum_{s=1}^\rho
   |k_{s-1}|^{\alpha r}|k_s|^{(1-2\alpha)r}\Big)^{\frac{1}{r}}
 \lesssim \frac{Q^2}{\rho},
 \end{align}
 where we used $$\sum_{l=1}^{s-1}|k_l|^\alpha \lesssim |k_{s-1}|^{\alpha}, \quad \Big(
   \sum_{s=1}^\rho
   (|k_{s-1}|^{\alpha r}|k_s|^{(1-2\alpha)r}\Big)^{\frac{1}{r}}\lesssim1.$$

\noindent Following the similar arguments as in \eqref{DLY-4.26} and
\eqref{DLY-4.29E}, we have
 $$\|\frac{Q^2}{\rho}\sum_{s=1}^\rho\!
 \int_0^t\!|k_s|^\alpha(\sum_{l\le s-1}\!|k_l|^\alpha
 S_{t-\tau}\mathbb{P}H_{s,l})d\tau\|_{L^2_TL^{q_\alpha}_x}\lesssim\frac{Q^2}{\rho}.
 $$

\noindent Estimates for the second part of \eqref{DLY-4.28}
 follow in the similar way, hence we obtain the desired results.
\end{proof}

\smallskip
\begin{remark}\label{D-1} {\rm
 By checking the proof of Lemmas \ref{Lem:3.9} and \ref{Lem:3.10},
 we obtain that the estimates for $u_{2,1}$ and $u_{2,2}$ are also
 true for any other $2\le q\le\infty$.
}
\end{remark}

\medskip
 Finally, we prove the lower bound of $u_{2,0}$ in critical space
 $\dot{F}^{-\alpha,r>2}_{q_\alpha}$ and the upper bound of $u_{2,0}$ in its well-posed space. Specially, the lower bound obtained plays
 a {\it crucial} role in the proof of norm inflation.

To obtain such bounds,  we will use several Fourier analysis methods. Due to the vector-valued nature of velocity field and the divergence free condition, we not only need to explore each of the three components but also need to analyze the action of Leray projection operator
 $\mathbb{P}$. 
 \medskip
 \begin{lemma}\label{Lem:3.11}
 For ${|k_0|^{-2\alpha}}\ll T\ll 1$, $q_\alpha=\frac{3}{\alpha-1}$, $1<\alpha<\frac{5}{4}$ and $2\le r\le \infty$, we
 get
 \begin{align}\label{DLY-4.34}
 &\|u_{2,0}(T)\|_{\dot{F}^{-\alpha,r}_{q_\alpha}}\gtrsim{Q}^{2},\\
 &\|u_{2,0}\|_{L^{q_\alpha}_xL^2_T}\lesssim
 T^{\frac{1}{2}}Q^2.\label{DLY-4.35}
 \end{align}
 \end{lemma}
\begin{proof}
 We first prove \eqref{DLY-4.34} and divide the proof into three steps.

 \noindent {\bf{Step 1}.}  From (H1)--(H2), \eqref{DLY-4.19} and \eqref{DLY-4A}, we have
 \begin{align}\label{DLY-4.36}
 \textrm{supp}\,\widehat{u}_{2,0}\!&\subset
 B_{\frac{1}{2}}(0)\cup{B}_{\frac{1}{2}}(k_s\!+\!k_s')
 \cup{B}_{\frac{1}{2}}(-(k_s\!+\!k_s')\,)
 \subset
 \Big\{\,|\xi|<\!\frac{8}{3}\,2^{2}\,\Big\},
 \end{align}
where $\pm(k_s+k_s')=(\pm7,0,0)$. Hence {\it for any $t\!>\!0$,
$u_{2,0}(x,t)\!\in\!\mathcal{C}^2(\mathbb{R}^3)$}.

 \noindent \textbf{Step 2}.
From $\dot{F}^{-\alpha,r}_{q_\alpha}\hookrightarrow
\dot{B}^{-\alpha,\infty}_{q_\alpha}$ and Definition \ref{Def:1.1},
we
 have
  \begin{align}
     \|u_{2,\,0}\|_{\dot{F}^{-\alpha,r}_{q_\alpha}}&\gtrsim\|u_{2,\,0}\|_{\dot{B}^{-\alpha,\infty}_{q_\alpha}}
     \gtrsim
       \|\Delta_2u_{2,0}\|_{L^{q_\alpha}_x}+ \|\Delta_3u_{2,0}\|_{L^{q_\alpha}_x} \nonumber\\
     &\gtrsim
 \|(\Delta_2\!+\!\Delta_3)u_{2,0}\|_{L^{q_\alpha}_x}\gtrsim
 \|(\Delta_2\!+\!\Delta_3)u_{2,0}\|_{L^{\infty}_x},\label{4.36-A}
  \end{align}
  where in the last inequality we used Bernstein's inequality.

 \noindent \textbf{Step 3}. In this step, it suffices to prove 
  \begin{align}\label{4.36-B}
  \Big((\Delta_2\!+\!\Delta_3)u_{2,0}^{[3]}\Big)(-{\pi}/{14},0,0,t)\gtrsim Q^2,
  \end{align}
  where $u_{2,0}^{[3]}$ denotes the third component of $u_{2,0}$.
  Once we prove \eqref{4.36-B}, then combining $u_{2,0}\in\mathcal{C}^2(\mathbb{R}^3)$ with
 \eqref{4.36-A}, we obtain that
  $$\|u_{2,0}\|_{\dot{F}^{-\alpha,r}_{q_\alpha}}\gtrsim \|(\Delta_2\!+\!\Delta_3)u_{2,0}^{[3]}\|_{L^\infty_x}\gtrsim Q^2,$$
  which is the desired \eqref{DLY-4.34}.

  To prove \eqref{4.36-B}, we recall that from Remark 3.7,   $\widehat{F}_{s}$ is purely imaginary, $\widehat{\mathbb{P}}$ and $u_{2,0}$
   are real-valued, $(\Delta_2\!+\!\Delta_3)u_{2,0}^{[3]}(x,t)$ can be rewritten as
 \begin{align}
  \label{DLY-4.37}
 (\Delta_2\!+\!\Delta_3)u_{2,0}^{[3]}(x,t)
 & =\frac{Q^2}{{(2\pi)^{3}}\rho}\sum_{s=1}^\rho\int_0^t\!\int_{\mathbb{R}^3}
   \frac{i\sin{(x\!\cdot\!\xi)}|k_s|^{2\alpha}}
   {e^{(t-\tau)|\xi|^{2\alpha}}}\,
   \widehat{\mathbb{P}}_3(\xi)\!\cdot\!(( \Delta_2 \!+\! \Delta_3 ) F_s)^\wedge (\xi) d\tau{d}\xi\nonumber\\
 :\!&=\frac{Q^2}{(2\pi)^3\rho}\sum_{s=1}^\rho I_{s}{(x,t)}
 \end{align}
where
$\widehat{\mathbb{P}}_3(\xi)=({\frac{-\xi_3\xi_1}{|\xi|^2},\frac{-\xi_3\xi_2}{|\xi|^2},\frac{|\xi_1|^2+|\xi_2|^2}{|\xi|^2}})$
 is the third row vector of $\widehat{\mathbb{P}}=I_d-\frac{\xi\otimes\xi}{|\xi|^2}$. 

 Later on, we will show that given $x_0=(-\frac{\pi}{14},0,0)$, there exists positive constant $\delta$ such that
 for any $1\le s\le \rho$ and $\frac{1}{|k_0|^{2\alpha}}\ll t\ll
 1$,
  $$I_{s}(x_0,t)\ge \delta .$$
   Noticing that $B_{\frac{1}{2}}(0)\cap\{\frac{3}{4}\ 2^2<|\xi|<\frac{8}{3}\ 2^3\}=\emptyset$ and
 \begin{align*}
 (\Delta_2\!+\!\Delta_3)F_s
 &= \!{S_\tau}{\Phi}_{k_s'} \!\cdot \!\nabla{S_\tau}{\Psi}_{k_s}
    \!+\!{S_\tau}\bar{\Phi}_{k_s'} \!\cdot\! \nabla{S_\tau}\bar{\Psi}_{k_s}
     \!+\!{S_\tau}{\Psi}_{k_s} \!\cdot \!\nabla{S_\tau}{\Phi}_{k_s'}
     \!+\!{S_\tau}\bar{\Psi}_{k_s} \!\cdot\!
      \nabla{S_\tau}\bar{\Phi}_{k_s'},
      \end{align*}
 then we get
  \begin{align}\label{Is}
I_{s}(x,t)& = \sum_{\ell=1}^4
     \int_{\mathbb{R}^3}\!
     \int_{\mathbb{R}^3}A_{s}(x,t,\xi,\eta)B_{\ell,\,s}(\xi,\eta){d}\eta{d\xi}
     :=\sum_{\ell=1}^4 I_{s,\ell}(x,t),
   \end{align}
where
 \begin{align}\label{As}
  A_{s}(x,t,\xi,\eta)&:=\sin(-{x\!\cdot\!\xi})
    \int_{0}^t\!
    {\,\,|k_s|^{2\alpha}} { e^{-t|\xi|^{2\alpha}+\tau(|\xi|^{2\alpha}-|\xi-\eta|^{2\alpha}-|\eta|^{2\alpha})}}
    d\tau\nonumber\\
    &\;=\frac{|k_s|^{2\alpha}\,\sin(-{x\!\cdot\!\xi})\, (e^{-t|\xi|^{2\alpha}}\!-\!e^{-t(|\xi-\eta|^{2\alpha}+|\eta|^{2\alpha})})}
    {|\xi\!-\!\eta|^{2\alpha}\!+\!|\eta|^{2\alpha}\!-\!|\xi|^{2\alpha}}
    \end{align}
    and
   \begin{align}
   \left\{\begin{aligned}
 &B_{1,s}(\xi,\eta)=  \widehat{\Phi}_{k_s'}(\xi\!-\!\eta)\cdot\eta\;\widehat{\mathbb{P}}_3(\xi)\!\cdot\!\widehat{\Psi}_{k_s}(\eta),\\
 &\hskip1.53cm=
    \widehat{\psi}_{k_s'}(\xi\!-\!\eta)\widehat{\psi}_{k_s}(\eta)
    \frac{(\xi_1\eta_2\!-\!\xi_2\eta_1)\,(\xi_1^2\eta_2\!+\!\xi_2^2\eta_2\!+\!\xi_3\xi_2 \eta_3)}{|\xi|^2|k_s|^{2}},\\
 &B_{3,s}(\xi,\eta) =   \widehat{\Psi}_{k_s}(\xi\!-\!\eta)\cdot\eta\;\widehat{\mathbb{P}}_3(\xi)\!\cdot\!\widehat{\Phi}_{k_s'}(\eta),\\
 &\hskip1.53cm=
    \widehat{{\psi}}_{k_s}(\xi\!-\!\eta)
    \widehat{{\psi}}_{k_s'}(\eta)
    \frac{ (\xi_3\eta_2-\xi_2\eta_3)\,(\xi_3\xi_2\eta_1-\xi_3\xi_1\eta_2)}{|\xi|^2|k_s|^{2}} ,\\
 &B_{2,s}(\xi,\eta)= -B_{1,s}(-\xi,-\eta),\;\;\;\;\;
  B_{4,s}(\xi,\eta)= -B_{3,s}(-\xi,-\eta).
  \end{aligned}\right.
  \label{Bs}
     \end{align}\
 From  \eqref{Is}--\eqref{Bs}, we observe that
 \begin{align}\label{Is12}
 I_{s,1}(x,t)=I_{s,2}(x,t),\quad I_{s,3}(x,t)=I_{s,4}(x,t).
 \end{align}

 Denoting $\tilde{\eta}:=\eta-k_s\in B_{{1}/{4}}(0)$, $\xi_0:=k_s+k_s'=(7,0,0)$ and recalling that
 $x_0=(-\frac{\pi}{14},0,0)$, we get
  \begin{align*}
   I_{s,1}(x_0,t)=&\int_{\mathbb{R}^3}\int_{\mathbb{R}^3}
   \widehat{\psi}(\xi-\!\xi_0-\!\tilde{\eta})\widehat{\psi}\,(\tilde{\eta})C_{s,1}(\xi,\tilde{\eta})\,A_{s}(x_0,t,\xi,\tilde{\eta}+k_s)d\xi d\tilde{\eta},
   \end{align*}
   where
   \begin{align}\label{cs}
   C_{s,1}(\xi,\tilde{\eta})=\frac{(\xi_1(\tilde{\eta}_2+|k_s|)\!-\!\xi_2\tilde{\eta}_1)\,((\xi_1^2+\xi_2^2)(\tilde{\eta}_2+|k_s|)\!+\!\xi_3\xi_2
   \tilde{\eta}_3)}{|\xi|^2|k_s|^{2}}.
  \end{align}

   Since $|\tilde{\eta}|<\frac{1}{4}$,
   $|\xi|<|\xi_0|\!+\!\frac{1}{2}$ and $2^{m_0}=|k_0|<|k_s|=2^{\frac{(s+1)(s+2m_0)}{2}}$, as $m_0$ tends to
   infinity, we have 
   $$\lim_{m_0\rightarrow\infty}C_{s,1}(\xi,\tilde{\eta})=\frac{\xi_1(\xi_1^2+\xi_2^2)}{|\xi|^2}$$
 and for $\frac{1}{|k_0|^{2\alpha}}\ll t\ll 1$,
  $$\lim_{m_0\rightarrow \infty}A_{s}(x_0,t,\xi,\tilde{\eta}+k_s)=\frac{e^{-t|\xi|^{2\alpha}}\sin(\frac{\pi \xi_1}{14})}{2}.$$
   As a consequence,
  \begin{align}\nonumber
  \lim_{m_0\rightarrow\infty} (I_{s,1} + I_{s,2})=
  \int_{\mathbb{R}^3}\int_{\mathbb{R}^3}
   \widehat{\psi}(\xi-\!\xi_0-\!\tilde{\eta})\widehat{\psi}\,(\tilde{\eta})
   \frac{\xi_1(\xi_1^2+\xi_2^2)\sin(\frac{\pi \xi_1}{14})}{|\xi|^{2}}e^{-t|\xi|^{2\alpha}}d\xi d\tilde{\eta}.
  \end{align}
   Similarly, we have
  \begin{align}
  \nonumber
  \lim_{m_0\rightarrow\infty} (I_{s,3} + I_{s,4})=
  \int_{\mathbb{R}^3}\int_{\mathbb{R}^3}
   \widehat{\psi}(\xi-\!\xi_0-\!\tilde{\eta})\widehat{\psi}\,(\tilde{\eta})
   \frac{-  \xi_1\xi_3^2\sin(\frac{\pi \xi_1}{14})}{|\xi|^{2}}e^{-t|\xi|^{2\alpha}}d\xi d\tilde{\eta}.
  \end{align}

 For any $\xi\in B_{1/2}(\xi_0)\cup B_{1/2}(-\xi_0)$,\;
 $\max\{|\xi_1\pm 7|,\, |\xi_2|,\, |\xi_3|\}<1/2$,
 there exists absolute positive constants $\delta$ and $N_0$ such
 that if $m_0>N_0$, then for any $1\le s\le \rho$,
   $$I_{s}(x_0,t)=\sum_{\ell=1}^4 I_{s,\ell}(x_0,t)\ge  \delta \int_{\mathbb{R}^3}\int_{\mathbb{R}^3}
   \widehat{\psi}(\xi-\!\xi_0-\!\tilde{\eta})\widehat{\psi}\,(\tilde{\eta})d\xi d\tilde{\eta}=\delta,$$
 which concludes the key estimate \eqref{4.36-B} by using
 \eqref{DLY-4.37}.

\medskip
 It remains to prove \eqref{DLY-4.35}, i.e. $\|u_{2,0}\|_{L^{q_\alpha}_xL^2_T}\lesssim\|u_{2,0}\|_{L^2_TL^{q_\alpha}_x}\lesssim
  Q^2T^{\frac{1}{2}}$. We need to estimate $u_{2,0}^{[1]}$, $u_{2,0}^{[2]}$ and $u_{2,0}^{[3]}$.
 Similar to \eqref{DLY-4.37}, we have
 \begin{align}\label{DLY-4.38A}
  \widehat{u}_{2,0}^{[3]}(\xi,t)&=
  \frac{Q^2}{{(2\pi)^{\frac{3}{2}}}\rho}\sum_{s=1}^\rho\int_0^t
    |k_s|^{2\alpha} {e^{-(t-\tau)|\xi|^{2\alpha}}}\,
   \widehat{\mathbb{P}}_3(\xi)\!\cdot\! \widehat{F}_s(\xi) d\tau.
  \end{align}
 By using divergence free condition of $\Psi_{k_s}$, $\Phi_{k_s'}$ and $|\xi|<\frac{8}{3}\,2^2$ as well as $\widehat{\psi}(\cdot)\ge0$,
 it is easy to show that, for instance $S_\tau\Psi_{k_s}\!\cdot\!\nabla S_\tau\Phi_{k_s'}=\nabla\!\cdot\!(S_\tau\Psi_{k_s}\otimes
 S_\tau\Phi_{k_s'})$.  Thus we get
 \begin{align}
 \label{DLY-4.38B}
 {e^{-(t-\tau)|\xi|^{2\alpha}}}|\widehat{\mathbb{P}}_3(\xi)\!\cdot\! \widehat{F}_s(\xi)|\lesssim
 e^{-c\tau|k_s|^{2\alpha}}
 f(\xi;k_s,k_s'),
 \end{align}
 where $
f(\xi;k_s,k_s'):=(\widehat{\psi}_{k_s}\ast\widehat{\psi}_{k_s'}+\widehat{\psi}_{k_s}\ast\widehat{\bar{\psi}}_{k_s}+
 \widehat{\psi}_{k_s'}\ast\widehat{\bar{\psi}}_{k_s'})(\xi).
$

  Combining \eqref{DLY-4.38A}--\eqref{DLY-4.38B}, we observe that
  for any $t>0$,
 \begin{align}\nonumber
 |\widehat{u}_{2,0}^{[3]}(\xi,t)| \lesssim \frac{Q^2}{\rho} \sum_{s=1}^\rho
  \int_{0}^t |k_s|^{2\alpha} e^{- c \tau |k_s|^{2\alpha} } d\tau f(\xi;k_s,k_s')\lesssim Q^2f(\xi;k_s,k_s').
 \end{align}
 Similarly, we get $|\widehat{u}_{2,0}^{[1]}|+|\widehat{u}_{2,0}^{[2]}|\lesssim Q^2 f(\xi;k_s,k_s')$. By applying
  Hausdorff-Young's inequality and  $\|\widehat{\psi}\|_{L^1_\xi\cap L^{\frac{q_\alpha}{q_\alpha-1}}_\xi}\lesssim1$
  in (H2) to $f(\xi;k_s,k_s')$, we have
 \begin{align*}
 \|u_{2,0}\|_{{L}^{2}_TL^{q_\alpha}_x}\lesssim
 \|\widehat{u}_{2,0}\|_{L^2_TL^{\!\frac{q_\alpha}{q_\alpha-1}}_\xi}\lesssim
   Q^2\Big\|\|f(\xi;k_s,k_s')\|_{L^{\frac{q_\alpha}{q_\alpha-1}}_\xi}\Big\|_{L^2_T}
  \lesssim  Q^2T^{\frac{1}{2}}.
  \end{align*}

 Therefore, we complete the proof.
  \end{proof}

\subsection{Estimates of remainder $y$}
 In this subsection, we use iteration arguments to prove the {\it a-priori} estimate for remainder $y$.
 Recall that $y$  satisfy the integral equations
 \eqref{eq:4.2}, i.e.
 $$y=G_0+G_1-G_2$$
  with initial condition $y|_{t=0}=0$ and
  \begin{align*}
  G_2 = B(y,y), \,
  G_1  = B(y,u_2\!-\!u_1)+B(u_2\!-\!u_1,y), \,
  G_0  = B(u_{2}, u_1\!-\!u_2)+B(u_1,u_{2}).
      \end{align*}

 From Lemma \ref{Lem:3.6}, we observe that in order to obtain more accurate decay estimate for $y$,
 it suffices to split $u_1$, $u_2$ and $u_{2}$ into
 two terms, e.g.
 \begin{align*}
 \left\{\begin{aligned}&u_1=u_1\chi_{_{[0,T_\sigma]} }(t) +
 u_1\chi_{_{[T_\sigma,T_{\sigma+1}]} }(t),\\
 &u_2=u_2\chi_{_{[0,T_\sigma]} }(t) +
  u_2\chi_{_{[T_\sigma,T_{\sigma+1}]} }(t),\\
 &y=y\chi_{_{[0,T_\sigma]} }(t) + y\chi_{_{[T_\sigma,T_{\sigma+1}]} }(t),
\end{aligned}\right.
\end{align*}
 Plugging the above decompositions of $u_1$, $u_2$ and $y$ into
 $G_0$, $G_1$ and $G_2$, we have the following iteration rules which play an important role in controlling $y$.
\begin{lemma}\label{Lem:3.12}
 If $y$ solves system \eqref{eq:4.2}--\eqref{eq:4.3}, then for any $\sigma=0,1,2,\cdots,Q^3$ and for large
 enough $\rho$ and $|k_0|$ we have
 \begin{align}\label{DLY:3.71}
 \|y\|_{X_{T_{\sigma+1}}}%
    \lesssim Q^{\sigma+3} (\, {\rho}^{-1}+ {|k_0|^{-\alpha}}\, ).
    \end{align}
 Moreover, for any $T>|k_0|^{-2\alpha}$, we have
 \begin{align}\label{DLY:3.72}
 \|y\|_{X_{T}}
     \lesssim Q^{3} ( {\rho}^{-1} + T^{\frac{1}{2}} )+Q^{Q^3+3} (\, {\rho}^{-1} +  {|k_0|^{-\alpha}}\, ).
    \end{align}
 \end{lemma}

 \begin{proof}
 Applying Lemma \ref{Lem:2.3} to \eqref{eq:4.2}--\eqref{eq:4.3}, we
 have the following bilinear estimates:
\begin{align}
   \|y\|_{X_{T_{\sigma+1}}}
   &\lesssim  \|u_2\|_{X_{T_{\sigma+1}}}(\|u_1\|_{X_{T_{\sigma+1}}} \!+\! \|u_2\|_{X_{T_{\sigma+1}}}) + (\|u_1\|_{X_{T_{\sigma+1}}}\!\!\!+\!\|u_2\|_{X_{T_{\sigma+1}}})\|y\|_{X_{T_\sigma}}\nonumber\\
   &\;\;\;+(\|u_1\|_{X_{[T_\sigma,T_{\sigma+1}]}}\!+\!\|u_2\|_{X_{[T_\sigma,T_{\sigma+1}]}})\|y\|_{X_{T_{\sigma+1}}}
   +\|y\|_{X_{T_{\sigma+1}}}^2.\label{DLY:3.73}
     \end{align}
 Recalling that for any $1\le\sigma\le\beta$, $T_\sigma\le T_\beta$. Then from Lemmas \ref{Lem:3.4}--\ref{Lem:3.11}, we get
 \begin{align}\label{DLY:3.74}
 \|u_2\|_{X_{T_{\sigma+1}}}\!\lesssim\!
 {Q^2}{{\rho}}^{-1} \!+ Q^2T^{\frac{1}{2}}_\beta,\;\;
 \|u_1\|_{X_{T_\sigma}} \!\lesssim\! Q,\;\; \|u_1\|_{X_{[T_\sigma,T_{\sigma+1}]}}
 \!\lesssim\! Q^{-\frac{1}{2}}.
 \end{align}
 Plugging \eqref{DLY:3.74} in \eqref{DLY:3.73}, and assuming that $\rho>Q^5$, $T_\beta=|k_0|^{-2\alpha}<Q^{-5}$, we have
 \begin{align}\label{DLY:3.75}
 \|y\|_{X_{T_{\sigma+1}}}
  &\!\lesssim ({Q^2}{{\rho}}^{-1}\!+\!Q^2T^{\frac{1}{2}}_\beta)(Q+{Q^2}{{\rho}}^{-1} \!+\! Q^2T^{\frac{1}{2}}_\beta)
             +(Q+{Q^2}{{\rho}}^{-1} \!+\! Q^2T^{\frac{1}{2}}_\beta)\|y\|_{X_{T_\sigma}}\nonumber\\
  &\;\;\;\;+(Q^{-\frac{1}{2}}+{Q^2}{{\rho}}^{-1}\!+Q^2T^{\frac{1}{2}}_\beta)\|y\|_{X_{T_{\sigma+1}}} + \|y\|_{X_{T_{\sigma+1}}}^2\nonumber\\
  &\lesssim Q^3(\, {{\rho}^{-1}} + {|k_0|^{-\alpha}}\,)
             +Q\|y\|_{X_{T_\sigma}} + Q^{-\frac{1}{2}}\|y\|_{X_{T_{\sigma+1}}} + \|y\|_{X_{T_{\sigma+1}}}^2.
 \end{align}
 Similarly, when $T>T_\beta$,
 by splitting $[0,T]$ into $[0,T_\beta]$ and $[T_\beta,T]$, then
 using
  Corollary \ref{Cor:3.7} and \eqref{DLY:3.73}--\eqref{DLY:3.75}, we get $\|u_1\|_{X_{[T_\beta,T]}}\lesssim
 Q\rho^{-\frac{1}{2}}$ and
 \begin{align}\label{DLY:3.76}
 \|y\|_{X_{T}}
 &\lesssim Q^3(\,{{\rho}^{-1}} +T^{\frac{1}{2}})
             +Q\|y\|_{X_{T_\beta}} + Q^{-\frac{1}{2}}\|y\|_{X_{T}} + \|y\|_{X_{T}}^2.
 \end{align}
 Lemma \ref{Lem:3.4} ensures that $\|y\|_{X_{T_0}}$ can be small
 enough
 since $T_0=|k_\rho|^{-2\alpha}$ and
\begin{align*}\|u_1\|_{X_{T_0}}\lesssim
  {Q}{{\rho}^{-\frac{1}{2}}}T^\frac{1}{2}_0|k_\rho|^\alpha\lesssim  {Q}{ {\rho}^{-\frac{1}{2}}}\end{align*} and $\rho$ is large
 enough. Thus iteration argument can be applied to \eqref{DLY:3.75}--\eqref{DLY:3.76} and hence we obtain the desired results.
 \end{proof}

 Making use of Lemmas \ref{Lem:2.3}
 and \ref{Lem:3.12}, we obtain the
 following estimate.
\begin{corollary}\label{Cor:3.13}
 For any  $r>2$, sufficiently large $\rho$ and $|k_0|$ such that $\rho\gg Q^{Q^3+2}$, $|k_0|^{-\alpha}\ll Q^{-Q^3-2}$
 and $|k_0|^{-2\alpha}<T\ll Q^{-4}$, we have
 \begin{align}\label{DLY:3.77}
   & \|y(T)\|_{\dot{F}^{-\alpha,r}_{q_\alpha}}
\ll Q^2.
 \end{align}
\end{corollary}
\begin{proof}
 From \eqref{eq:4.2}--\eqref{eq:4.3}, we notice that $y(T)=G_0(T)+G_1(T)-G_2(T)$
 and $G_i(T)$ are several bilinear terms. By applying Lemma
 \ref{Lem:2.3}, we obtain that
 \begin{align}
 \|y(T)\|_{\dot{F}^{-\alpha,r}_{q_\alpha}}
   &\lesssim \|y(T)\|_{\dot{F}^{-\alpha,2}_{q_\alpha}}\lesssim   \|y\|_{L^\infty_T\dot{F}^{-\alpha,2}_{q_\alpha}}\nonumber\\
   &\lesssim
   \|u_2\|_{X_T}(\|u_1\|_{X_T}+\|u_2\|_{X_T})+(\|u_1\|_{X_T}+\|u_2\|_{X_T})\|y\|_{X_T}+\|y\|_{X_T}^2.\nonumber
 \end{align}
  Applying Lemmas \ref{Lem:3.4},
 \ref{Lem:3.9}--\ref{Lem:3.12} to the above inequality, we have
 \begin{align}
  \|y(T)\|_{\dot{F}^{-\alpha,r}_{q_\alpha}}
  &\lesssim( {Q^2}{\rho}^{-1} +  Q^2T^{\frac{1}{2}})(Q +  {Q^2}{\rho}^{-1} + Q^2T^{\frac{1}{2}})
  \nonumber\\
  &\;\;\;\;+(Q +   {Q^2}{\rho}^{-1} + Q^2T^{\frac{1}{2}})
  \Big(Q^{3}( {\rho}^{-1} + T^{\frac{1}{2}})+Q^{Q^3+3}( {\rho}^{-1} + {|k_0|^{-\alpha}})\Big)\nonumber\\
  &\;\;\;\;+\Big(Q^{3}( {\rho}^{-1} +\!T^{\frac{1}{2}})+Q^{Q^3+3}(\, {\rho}^{-1} + {|k_0|^{-\alpha}})\Big)^2
  \nonumber\\
  &\ll Q^2.\nonumber \end{align}
 Hence we prove the desired result.
  \end{proof}

 \subsection{Proof of Theorem
 \ref{Thm:1.3}}
  In this subsection, combining the results proved in Subsections 3.1--3.4, we are ready to prove the ill-posedness of the gNS by
  showing norm inflation.

\medskip
 \noindent\textit{Proof of Theorem
 \ref{Thm:1.3}}. Combining the equalities \eqref{eq:4.1} and \eqref{DLY-4.19}--\eqref{DLY-4.21}, the estimates \eqref{DLY-4.10},
 \eqref{DLY-4.22}, \eqref{DLY-4.27},  \eqref{DLY-4.34} and \eqref{DLY:3.77}, we have
 \begin{align}
   \|u(T)\|_{\dot{F}^{-\alpha,r}_{q_\alpha}}&\ge
   \|u_{2,0}(T)\|_{\dot{F}^{-\alpha,r}_{q_\alpha}}\nonumber\\
   &\;\;-\Big(\|u_1(T)\|_{\dot{F}^{-\alpha,r}_{q_\alpha}}\!+\!\|u_{2,1}(T)\|_{\dot{F}^{-\alpha,r}_{q_\alpha}}\!\!+\!
   \|u_{2,2}(T)\|_{\dot{F}^{-\alpha,r}_{q_\alpha}}\!+\!\|y(T)\|_{\dot{F}^{-\alpha,r}_{q_\alpha}}\Big)\nonumber\\
   &\gtrsim
   \|(\Delta_2+\Delta_3)u_{2,0}^{[3]}(T)\|_{L^{\infty}_x}\nonumber\\
   &\;\;-\Big(\|u_1(T)\|_{\dot{F}^{-\alpha,r}_{q_\alpha}}\!+\!\|u_{2,1}(T)\|_{\dot{F}^{-\alpha,2}_{q_\alpha}}\!
   +\!\|u_{2,2}(T)\|_{\dot{F}^{-\alpha,2}_{q_\alpha}}\!+\!\|y(T)\|_{\dot{F}^{-\alpha,2}_{q_\alpha}}\Big)\nonumber\\
   &\gtrsim Q^2\Big(1-{Q}^{-1}\rho^{\frac{1}{r}-\frac{1}{2}}-\rho^{-1}-o(1)\Big)
   \gtrsim Q^2,\nonumber\label{DLY:3.78}
 \end{align}
 where $0<o(1)\ll\frac{1}{2}$, $\rho\!\gg\!Q^{Q^3+2}$ and $|k_0|^{-2\alpha}\!<\!T\!\ll \!Q^{-4}$.
  Hence we finish the proof.

\vspace*{3ex}
 \noindent\textbf{\Large Appendix}

\vspace*{1ex}


In this appendix, we will give a proof of Lemma 2.4 and state some
extensions for reader's convenience. In fact, the equivalent
estimates (\ref{DLY:2.6}) of  Lemma 2.4 can be immediately concluded
from the following general Littlewood-Paley g-function
characterizations of $L^p(\mathbb{R}^n)$, which in turn base on the
vector-value singular integrals theory, see e.g. Stein \cite[p.46
and p.185]{Stein:1993}.

\noindent\textbf{Lemma A.1.}\,
 \label{Lem:A.1}
 {\it Let $\Phi(x)$ be any  function on $\mathbb{R}^n$ satisfying $\int_{\mathbb{R}^n}\Phi(x)dx=0$ and
\begin{equation}\tag{A.1}|\Phi(x)|+|\nabla\Phi(x)|\le A(1+|x|)^{-n-1},\end{equation}
for some constant $A>0$.  Then for any $1<q<\infty$ the following
estimate
\begin{equation}\tag{A.2}
\|s_\Phi f\|_{L^q_x}:=\Big\|\Big(\int_0^\infty |\Phi_t\ast
f|^2\frac{dt}{t}\Big)^{1/2}\Big\|_{L^q_x}\le
C_q\|f\|_{L^q}
\end{equation} holds for
$\Phi_t(x)=\frac{1}{t^{n}}\Phi(\frac{x}{t})$  with $t>0$.
Furthermore, if $\Phi$ is nondegenerate, in the sense that there
exists a function $\Psi$ satisfying the same conditions as $\Phi$
such that
\begin{equation}\tag{A.3}
\int_0^\infty
\widehat{\Phi}(t\xi)\widehat{\Psi}(t\xi)\frac{dt}{t}=1,\ \  \xi\neq
0,
\end{equation}
then the converse inequality $\|f\|_{L^q_x}\le C'_q\|s_\Phi
f\|_{L^q_x}$ holds for any $1<q<\infty$. }

\vskip 0.2cm \noindent\textbf{Proof of Lemma A.1}
 Let
$\Phi(x)=\mathcal{F}^{-1}(|\cdot|^\alpha e^{-|\cdot|^{2\alpha}})(x)$
for any $\alpha\ge1$. Then
$$\int_{\mathbb{R}^n}\Phi(x)dx=\widehat{\Phi}(0)=0$$
and it is easy to check that $\Phi(x)$ satisfies the condition
({A.1}) by the similar argument as done in Lemma 2.2. Moreover, we
can choose $\Psi=c \Phi$ for some $c>0$ such that
$$\int_0^\infty |\widehat{\Phi}(t\xi)|^2\frac{dt}{t}=1/c,\ \  \xi\neq 0,$$
which means the $\Phi$ is non-degenerate. Thus it follows from
Theorem \ref{Lem:A.1} that
\begin{equation}\tag{A.4}
\Big\|\Big(\int_0^\infty
|t^\alpha(-\Delta)^{\alpha/2}e^{-t^{2\alpha}(-\Delta)^{\alpha}}
h|^2\frac{dt}{t}\Big)^{1/2}\Big\|_{L^q_x}\sim\|h\|_{L^q_x}
\end{equation}
Equivalently, by the variable $t$ changing and set
$h=(-\Delta)^{-\alpha/2}f$,  we have
\begin{equation}\nonumber
\Big\|\Big(\int_0^\infty |e^{-t(-\Delta)^{\alpha}}f|^2
dt\Big)^{1/2}\Big\|_{L^q_x}\sim \|(-\Delta)^{-\alpha/2}f\|_{L^q_x}.
\end{equation}
Note that by Littlewood-Plaey theorem and isomorphism, it follows
that  $$\|(-\Delta)^{-\alpha/2}f\|_{L^q_x} \sim
\|(-\Delta)^{-\alpha/2}f\|_{\dot{F}^{0,2}_q}\sim
\|f\|_{\dot{F}^{-\alpha,2}_q}$$ for any $1<q<\infty$, hence, the
fractional semigroup characterization
 \begin{equation}\tag{A.5}
\Big\|\Big(\int_0^\infty
|e^{-t(-\Delta)^{\alpha}}f|^2dt\Big)^{1/2}\Big\|_{L^q_x}\sim
\|f\|_{\dot{F}^{-\alpha,2}_q}
\end{equation}
holds for any $1<q<\infty$.  When $n=3$, $1<\alpha<\frac{5}{4}$ and
$q_\alpha=\frac{3}{\alpha-1}$, we can immediately get the desired
results in Lemma 2.4.

\vskip0.2cm \begin{remark}\label{Rem:A.2} {\rm Intrinsically, we can
extend also the estimate (A.5) to general case,  for instance, for
any $s<0$ and $1<r, q<\infty$,
 \begin{equation}\tag{A.6}
\Big\|\Big(\int_0^\infty
|t^{-\frac{s}{2\alpha}}e^{-t(-\Delta)^{\alpha}}f|^r\frac{dt}{t}\Big)^{1/r}\Big\|_{L^q_x}\sim
\|f\|_{\dot{F}^{s,r}_q}.
\end{equation}
holds.  In particular, when $s=-\alpha$ and $r=2$,  we immediately
obtain the estimate (A.5) above. However, it should be pointed out
that the proof of general estimate (A.6) is different and more
involved than the special index $r=2$, essentially depending on a
vector-valued version of maximal functions inequality, originally
due to Fefferman and Stein. In this deep connection, one can see
Triebel book \cite[p. 101]{Triebel:1978} for many general
characterizations of nonhomogeneous Triebel-Lizorkin space
$F^{s,r}_q(\mathbb{R}^n)$, where  one can check similar methods also
work well for the proof of the homogeneous type (A.6). Hence we omit
these details in the appendix for concision.}
\end{remark}

 \vspace*{1ex}
 \noindent {\bf Acknowledgmens:}
 Chao Deng is supported by PAPD of Jiangsu Higher Education Institutions and JSNU (No.\,9212112101), and NSFC Tianyuan Found (No.\,11226180); He is also partially supported by the NSFC (No.\,11171357, 11271166);
 Part of this work is done when Chao is visiting Penn State University and he would like to express his gratitude to
 professor Chun Liu and the Math Department of PSU for their hospitality.  Xiaohua Yao
 is supported by NSFC (No.\,10801057), NCET-10-0431 and the Special
 Fund for Basic Scientific Research of Central Colleges (No.\,CCNU12C01001).
 
 \vspace*{1ex}
 {{\small \hspace*{5cm}Chao Deng \;\quad {\it{deng315@yahoo.com.cn},}\\
 {\small \hspace*{5cm}Department of Mathematics,}\\
 {\small \hspace*{5cm}Jiangsu Normal University, Xuzhou 221116, PRC}\\
   \\
{\small \hspace*{5cm}Xiaohua Yao\quad {\it{yaoxiaohua@mail.ccnu.edu.cn},}\\
 {\small \hspace*{5cm}Department of  Mathematics,}\\
 {\small \hspace*{5cm}Central China Normal University, Wuhan 430079, PRC}\\
}
 \date{}

\end{document}